\documentclass[11pt]{amsart}

\usepackage{amssymb,amsmath,amsxtra,amscd,eucal} 
\theoremstyle{plain}
\numberwithin{equation}{section}
\newtheorem{theorem}{Theorem}
\newtheorem{proposition}{Proposition}
\newtheorem{corollary}{Corollary}
\newtheorem{lemma}{Lemma}
\numberwithin{lemma}{section}
\numberwithin{theorem}{section}
\numberwithin{corollary}{section}
\numberwithin{proposition}{section}

\theoremstyle{definition}
\newtheorem{definition}{Definition}
\numberwithin{definition}{section}

\theoremstyle{remark}

\numberwithin{example}{section}
\newtheorem{remark}{Remark}

\newcommand{\nc}{\newcommand}
\renewcommand{\H}{\ensuremath{\mathcal{H}}}
\nc{\Hs}{\ensuremath{\mathcal{H}^{\sigma}}}
\nc{\X}{\ensuremath{X}}
\nc{\Y}{{\mc Y}}
\nc{\Z}{{\mathbb Z}}
\nc{\nul}{\ensuremath{\varnothing}}
\nc{\dcox}{h^\vee}

\nc{\sdp}{\ensuremath{\ltimes}}

\nc{\znz}[1]{\ensuremath{\mathbb{Z} / #1 \mathbb{Z}}}

\nc{\Xnot}{\overset{\circ}{X}}
\nc{\vbun}{\ensuremath{\mathcal{V}^{H}_{\Xnot}}}
\nc{\vbund}{\ensuremath{\mathcal{V}^{H,*}_{\Xnot}}}

\nc{\fps}[1]{\ensuremath{\mathbb{C}[[#1]]}}
\nc{\ls}[1]{\ensuremath{\mathbb{C}((#1))}}
\nc{\zf}[1]{\ensuremath{z^{\frac{1}{#1}}}}
\nc{\wf}[1]{\ensuremath{w^{\frac{1}{#1}}}}
\nc{\tf}[1]{\ensuremath{t^{\frac{1}{#1}}}}
\newcommand{\ru}[1]{\ensuremath{e^{\frac{2 \pi i}{N}}}}

\nc{\Orb}[1]{\ensuremath{\mathbf{O}_{#1}}}

\nc{\Cnot}{\ensuremath{C \backslash \{ \nu^{-1}(x_{i}) \}}}

\nc{\g}{\ensuremath{\mathfrak{g}}}
\nc{\goutc}{\ensuremath{\g^{\sigma}_{\on{out}}(\Caff)}}
\nc{\Lg}{\ensuremath{L \g}}
\nc{\ghat}{\ensuremath{\hat{\g}}}

\nc{\Ksto}{\ensuremath{\mathcal{K}^{*}_{odd}}}
\nc{\Osto}{\ensuremath{\mathcal{O}^{*}_{odd}}}
\nc{\Ost}{\ensuremath{\mathcal{O}^{*}}}

\nc{\Vkg}{\ensuremath{V_{k}(\mathfrak{g})}}

\nc{\vac}{|0\rangle}
\nc{\V}{\ensuremath{\mathcal{V}}}
\nc{\Vt}{\ensuremath{\wt{\mathcal{V}}}}

\nc{\nV}{\nu^*(\V)}

\nc{\Vtx}{\ensuremath{\widetilde{\mathcal{V}}_{X}}}
\nc{\Ptx}{\ensuremath{\Pi_{C}}}
\nc{\Mt}{\ensuremath{M^{\sigma}}}
\nc{\Mtp}[1]{\ensuremath{\mathcal{M}^{\sigma}_{#1}}}
\nc{\Mp}[1]{\ensuremath{\mathcal{M}_{#1}}}
\nc{\Np}[1]{\ensuremath{\mathcal{N}_{#1}}}
\nc{\Mtc}{\ensuremath{\mathcal{M}}}
\nc{\Ym}[1]{\ensuremath{Y^{#1}}}

\nc{\fd}[1]{\ensuremath{\mathcal{D}_{#1}}}
\nc{\fpd}[1]{\ensuremath{\mathcal{D}^{\times}_{#1}}}
\nc{\AutO}{\ensuremath{\Aut \mathcal{O}}}
\nc{\AutOr}[1]{\Aut_{#1} \mathcal{O}}

\nc{\AutD}[1]{\Au(\mathcal{D}_{#1})}
\nc{\AutDN}[2]{\Au_{#1}(\mathcal{D}_{#2})}
\nc{\AutX}[1]{\ensuremath{\emph{\mc{A}ut}_{#1}}}

\nc{\discs}[1]{\ensuremath{\coprod_{p \in \nu^{-1}(x_{#1})} \fpd{p}}}

\newcommand{\Cnotp}

\newcommand{\Houtc}{ \ensuremath{ \mathcal{H}_{\on{out}}(\Caff) }}
\newcommand{\Houtcp}{ \ensuremath{ \mathcal{H}_{\on{out}}(\Caffp) }}

\nc{\Caffp}{\ensuremath{C^{'}_{aff}}}
\nc{\Caff}{\ensuremath{C_{aff}}}
\nc{\Csq}{\ensuremath{\hat{C}^{2}}}
\nc{\Otil}{\ensuremath{\widetilde{\mathcal{O}}}}
\nc{\Ctil}{\ensuremath{\widetilde{\mathcal{C}}}}
\nc{\phitil}{\ensuremath{\widetilde{\phi}}}

\nc{\pullPi}{\ensuremath{\pi^{*} \Pi_{X}}}

\nc{\Cb}[1]{\ensuremath{\wt{\mc{C}}_{\pi}(X,\{ x_{i} \}, 
    \pi_{x_{i}} )_{i=1 \cdots #1} }}

\nc{\Hco}[1]{\ensuremath{\wt{\mc{H}}_{\pi}(X,\{ x_{i} \}, 
    \pi_{x_{i}} )_{i=1 \cdots #1} }}

\nc{\Cbq}[1]{\ensuremath{\wt{\mc{C}}_{\pi}(X,\{ x_{i}, \nu(q) \}, 
    \pi_{x_{i}}, \pi_{\nu(q)} )_{i=1 \cdots #1} }}

\nc{\Coinv}[2]{\ensuremath{ \widetilde{\mathcal{H}}_{\pi}(C, \{ p_{i} \},
	\{ q_{j }\}, \{ \pi^{\sigma}_{i} \}, \{ \pi^{\lambda_{j}}_{j}
	\} )^{ i = 1 \cdots #1}_{j=1 \cdots #2} }}

\nc{\Cv}[1]{\ensuremath{\wt{\mc{\H}}_{\pi}(X,\{ x_{i} \}, \pi_{x_{i}}
    )_{i=1 \cdots m} }}

\nc{\points}[1]{\ensuremath{ \{x_{i} \}_{i=1 \cdots #1} }}

\nc{\Cvb}[2]{\ensuremath{\mc{C}_{#1}(X,\{x_{i} \}, \mc{M}_{x_{i}})_{i=
1 \cdots #2} }}

\nc{\Cvl}[2]{\ensuremath{\mc{C}^{\mc{L}}_{#1}(X,\{x_{i} \},
\mc{M}^{\mc{L}}_{x_{i}})_{i= 1 \cdots #2} }}

\nc{\Cl}[2]{\ensuremath{\mc{H}^{\mc{L}}_{#1}(X,\{x_{i} \},
\mc{M}^{\mc{L}}_{x_{i}})_{i= 1 \cdots #2} }}

\nc{\Coinvb}[2]{\ensuremath{\mc{H}_{#1}(X,\{x_{i} \}, \mc{M}_{x_{i}})_{i=
1 \cdots #2} }}

\nc{\on}{\operatorname}
\nc{\Spec}{\on{Spec}}
\nc{\Aut}{\ensuremath{\on{Aut}}}
\nc{\Der}{\on{Der}}
\nc{\End}{\on{End}}
\nc{\Au}{\ensuremath{ {\mc A}ut}}
\nc{\mc}{\mathcal}
\nc{\ol}[1]{\ensuremath{\overset{\circ}{#1}}}
\nc{\ep}{\epsilon}
\nc{\C}{{\mathbb C}}
\nc{\OO}{{\mc O}}
\nc{\K}{{\mc K}}
\nc{\DD}{{\mc D}}
\nc{\bD}{\overline{\DD}}
\nc{\pa}{\partial}
\nc{\secref}[1]{Section~{\ref{#1}}}
\nc{\corref}[1]{Corollary~{\ref{#1}}}
\nc{\wt}{\widetilde}
\nc{\Res}{\on{Res}}
\nc{\Ind}{\on{Ind}}
\nc{\nn}{_{(n)}}
\nc{\cA}{{\mc A}}
\nc{\Wick}{{\mathbf :}}
\nc{\Hom}{\on{Hom}}
\nc{\al}{\alpha}
\nc{\mb}{\mathbf}
\nc{\ds}{\displaystyle}

\begin{document}

\title{Twisted Modules over Vertex Algebras on Algebraic Curves}

\author{Edward Frenkel}\thanks{Partially supported by grants from the
Packard Foundation and NSF and by an NSERC graduate fellowship.}
\address{Department of Mathematics  
         University of California, Berkeley, CA 94720}
\email{frenkel@math.berkeley.edu}
\author{Matthew Szczesny}
\email{szczesny@math.berkeley.edu}
\date{December 2001; Revised October 2002}

\begin{abstract}

We extend the geometric approach to vertex algebras developed by the
first author to twisted modules, allowing us to treat orbifold models
in conformal field theory. Let $V$ be a vertex algebra, $H$ a finite
group of automorphisms of $V$, and $C$ an algebraic curve such that $H
\subset \on{Aut}(C)$. We show that a suitable collection of twisted
$V$--modules gives rise to a section of a certain sheaf on the
quotient $X=C/H$. We introduce the notion of conformal blocks for
twisted modules, and analyze them in the case of the Heisenberg and
affine Kac-Moody vertex algebras. We also give a chiral algebra
interpretation of twisted modules.
 
\end{abstract}

\maketitle

\section{Introduction}

Conformal field theory (CFT) in two dimensions provides a rich setting
in which several areas of mathematics such as representation theory
and algebraic geometry interact in a natural way. In recent years,
much effort has been spent on setting up a precise mathematical
framework for CFT. The algebraic aspect of the theory has been
formalized in the language of vertex algebras (see
\cite{Bo,FLM,KAC,FBZ}). In order to understand the rich geometry behind
CFT, this algebraic approach must be combined with a geometric
formalism.

In \cite{FBZ}, an algebro-geometric approach to vertex algebras is
introduced (see \cite{H,BD} for other approaches). Starting with a
conformal vertex algebra $V$ and an algebraic curve $\X$, one can
construct a vector bundle $\V_{X}$ on $X$ such that vertex operators
become (endomorphism-valued) sections of $\V^{*}_{X}$. This gives a
coordinate-free description of vertex operators and allows one to make
contact with the fascinating geometry pertaining to $\X$ and related
moduli spaces.

If a vertex algebra has a group of automorphisms, then its
representation theory may be enhanced by the inclusion of twisted
modules. The systematic study of twisted modules was initiated in
\cite{FLM} where twisted vertex operators were used in the
construction of the Moonshine Module vertex algebra (see Chapter 9 of
\cite{FLM} and the works \cite{Lep1,Lep2}). The notion of the twisted
module was formulated in \cite{FFR,Dong} following \cite{FLM}. Twisted
modules (or twisted sectors as they are known in the physics
literature) appear as important ingredients of the so-called orbifold
models of conformal field theory (see \cite{DHVW,DVVV}). They have
been extensively studied in recent years (see, e.g.,
\cite{LI,DLM,BKT}).

In this paper we extend the geometric formalism developed in
\cite{FBZ} to twisted modules over vertex algebras. Let $C$ be a
smooth projective curve, and $H \subset Aut(C)$ a finite group of
automorphisms of $C$ such that the stabilizer of the action of $H$ on
at a generic point of $C$ consists of the identity element of
$H$. Suppose furthermore that $V$ is a conformal vertex algebra, and
that $H$ acts on $V$ by conformal automorphisms. We show that with
these data, the vector bundle $\V_{C}$ acquires an $H$--equivariant
structure, lifting the action of $H$ on $C$. Let $X = C / H$ be the
quotient curve, and $\nu: C \to X$ the quotient map, ramified at the
fixed points of $H$. Denote by $\ol{C} \subset C$ the locus of points
in $C$ whose stabilizer in $H$ is the identity element. Let $\ol{X}
\subset X$ be the image of $\ol{C}$ in $X$ and $\ol{\nu}: \ol{C} \to
\ol{X}$ the restriction of $\nu$ to $\ol{C}$. Thus, $\ol{C}$ is a
principal $H$--bundle over $\ol{X}$.

The vector bundle $\V_{\ol{C}}$ over $\ol{C}$ carries an
$H$--equivariant structure and hence descends to a vector bundle
on $\ol{X}$ which we denote by $\V^{H}_{\ol{X}}$.
%Explicitly,
%\[
%\V^{H}_{ \ol{X} } = \Au_{\ol{C}} \underset{H \times \AutO}{\times} V  
%\]

Let $x \in X$. Then $x$ corresponds to an $H$--orbit ${\mb O}_{x}$ in
$C$. For each point $p \in \nu^{-1}(x)$, the stabilizer $H_{p}$ is a
cyclic group, which has a canonical generator $h_{p}$, the monodromy
around $p$ (generically, $H_p = \{ e \}$ and $h_p = e$). We call a
collection $\{ M^{h_{p}}_{p}\}$ of $h_{p}$--twisted modules
satisfying certain compatibilities, a $V$--module along $\nu^{-1}(x)$.
For example, if $h_{p}=e$, then each $M^{h_{p}}_{p}$ is an ordinary
$V$--module and the requirement is that if $p' = g(p) \in
\nu^{-1}(p)$, then $M^{h_{p'}}_{p'}$ is obtained from $M^{h_{p}}_{p}$
by twisting the $V$--action by the automorphism of $V$ corresponding
to $g$. If, on the other hand, $H = \Z/N\Z$ and $h_p$ is a generator
of $H$, then $M^{h_p}_p$ can be an arbitrary $h_p$--twisted
$V$--module.

We attach to a $V$--module ${\mc M}_x$ along $\nu^{-1}(x)$ a section
$\mc{Y}^{\mc{M}_{x}}$ of $\V^{H,*}_{\ol{X}}$ on $\fpd{x}$, the
punctured disc at $x$. Using this structure we define the spaces of
conformal blocks in the twisted setting. The space of conformal blocks
is associated to a pair $(C,H)$ as above and a collection of
$V$--modules along $\nu^{-1}(x)$ attached to a set of points of $X
\backslash \ol{X}$, and a (possibly empty) collection of $V$--modules
along $\nu^{-1}(x), x \in \ol{X}$.  We give two equivalent definitions
of the space of conformal blocks: using the action of a certain Lie
algebra obtained from Fourier coefficients of vertex operators, and
using analytic continuation (as in \cite{FBZ}).  In the case of the
Heisenberg and affine Kac-Moody vertex algebras this definition may be
simplified using twisted versions of the Heisenberg and affine Lie
algebras, respectively.

Finally, we explain the connection with the chiral algebra
formalism. The right ${\mc D}_X$--module $\cA = \V_X \otimes \Omega_X$
is a chiral algebra on $X$ in the sense of A. Beilinson and
V. Drinfeld \cite{BD} (see \cite{FBZ}, Ch. 18). The action of $H$ on
$V$ induces an action of $H$ by automorphisms of $\cA$. Then the twist
$\cA^{\ol{C}}$ of $\cA|_{\ol{X}}$ by the $H$--torsor $\ol{C}$,
$\cA^{\ol{C}} = \cA|_{\ol{X}} \underset{H}\times \ol{C}$, is also a
chiral algebra. Twisted $V$--modules correspond to
$\cA^{\ol{C}}$--modules supported at the points $x \in X \backslash
\ol{X}$, and the above space of conformal blocks may be defined in
terms of these $\cA^{\ol{C}}$--modules.

\medskip

\noindent {\bf Acknowledgments.} We thank D. Gaitsgory for a useful
discussion of the chiral algebra interpretation of twisted modules.

\section{Vertex algebras and modules}    \label{algebras and modules}

In this paper we will use the language of vertex algebras, their
modules, and twisted modules. For an introduction to vertex algebras
and their modules \cite{FLM,KAC,FBZ}, and for background on twisted
modules, see \cite{FFR,Dong,DLM}.

We recall that a conformal vertex algebra is a $\Z_+$--graded vector
space $$V = \bigoplus_{n=0}^\infty V_n,$$ together with a vacuum
vector $\vac \in V_0$, a translation operator $T$ of degree $1$, a
conformal vector $\omega \in V_2$ and a vertex operation
\begin{align*}
Y: V &\to \on{End} V[[z^{\pm 1}]], \\
A &\mapsto Y(A,z) = \sum_{n \in \Z} A_{(n)} z^{-n-1}.
\end{align*}
These data must satisfy certain axioms (see \cite{FLM,KAC,FBZ}). In
what follows we will denote the collection of such data simply by $V$.

A vector space $M$ is called a $V$--module if it is equipped with an
operation
\begin{align*}
Y^M: V &\to \on{End} M[[z^{\pm 1}]], \\
A &\mapsto Y^M(A,z) = \sum_{n \in \Z} A^M_{(n)} z^{-n-1}
\end{align*}
such that for any $v \in M$ we have $A^M_{(n)} v = 0$ for large enough
$n$. This operation must satisfy the following axioms:

\begin{itemize}

\item $Y^M(\vac,z) = \on{Id}_M$;

\item For any $v \in M$ there exists an element $$f_v \in
M[[z,w]][z^{-1},w^{-1},(z-w)^{-1}]$$ such that the formal power series
$$Y^M(A,z) Y^M(B,w) v \qquad \on{and} \qquad Y_M(Y(A,z-w) B,w) v$$ are
expansions of $f_v$ in $M((z))((w))$ and $M((w))((z-w))$,
respectively.

\end{itemize}

The power series $Y^M(A,z)$ are called vertex operators. We write the
vertex operator corresponding to $\omega$ as
\[
	Y^M(\omega,z) = \sum_{n \in \mathbb{Z}} L^M_{n} z^{-n-2},
\]
where $L^M_n$ are linear operators on $V$ generating the Virasoro
algebra. Following \cite{Dong}, we call $M$ \emph{admissible} if
$L^{M}_{0}$ acts semi-simply with integral eigenvalues.

Now let $\sigma_{V}$ be a conformal automorphism of $V$, i.e., an
automorphism of the underlying vector space preserving all of the
above structures (so in particular $\sigma_{V}(\omega) = \omega$). We
will assume that $\sigma_{V}$ has finite order $N>1$. A vector space
$M^\sigma$ is called a $\sigma_V$--{\em twisted} $V$--module (or
simply twisted module) if it is equipped with an operation
\begin{align*}
Y^{M^\sigma}: V &\to \on{End} M^\sigma[[z^{\pm \frac{1}{N}}]], \\ A
&\mapsto Y^{M^{\sigma}}(A,\zf{N}) = \sum_{n \in \frac{1}{N}\Z}
A^{M^\sigma}_{(n)} z^{-n-1}
\end{align*}
such that for any $v \in M^{\sigma}$ we have $A^{M^\sigma}_{(n)} v =
0$ for large enough $n$. Please note that we use the notation
$Y^{M^{\sigma}}(A,\zf{N})$ rather than $Y^{M^{\sigma}}(A,z)$ in the
twisted setting. This operation must satisfy the following axioms (see
\cite{FFR,Dong,DLM,LI}):

\begin{itemize}

\item $Y^{M^{\sigma}}(\vac,\zf{N}) = \on{Id}_{\Mt}$;

\item For any $v \in \Mt$, there exists an element 
\[
f_{v} \in \Mt [[\zf{N}, \wf{N} ]][z^{- \frac{1}{N}}, w^{-
\frac{1}{N}},(z-w)^{-1}]
\]
such that the formal power series $$Y^{\Mt}(A,\zf{N}) Y^{\Mt}(B,
\wf{N})v \qquad \on{and} \qquad Y^{\Mt}(Y(A,z-w) B, \wf{N})v$$ are
expansions of $f_{v}$ in $\Mt((\zf{N}))((\wf{N}))$ and
$\Mt((\wf{N}))((z-w))$, respectively.
%For any $v \in M$ there exists an element of $f_v \in
%M[[z,w]][z^{-1},w^{-1},(z-w)^{-1}]$ such that the formal power series
%$Y^M(A,z) Y^M(B,w) v$ and $Y_M(Y(A,z-w) B,w) v$ are expansions of
%$f_v$ in $M((z))((w))$ and $M((w))((z-w))$, respectively;

\item If $A \in V$ is such that $\sigma_V(A) = e^{\frac{2\pi i m}{N}}
A$, then $A^{M^\sigma}_{(n)} = 0$ unless $n \in \frac{m}{N} + \Z$.

\end{itemize}

The series $Y^{M^\sigma}(A,z)$ are called twisted vertex operators.
In particular, the Fourier coefficients of the twisted vertex operator
\[
	Y^{M^\sigma}(\omega,\zf{N}) = \sum_{n \in \mathbb{Z}}
	L^{M^\sigma}_{n} z^{-n-2},
\]
generate an action of the Virasoro algebra on $M^\sigma$. The
$\sigma_{V}$--twisted module $\Mt$ is called \emph{admissible} if
$L^{\Mt}_{0}$ acts semi-simply with eigenvalues in $\frac{1}{N}\Z$.

One shows in the same way as in \cite{FBZ}, Sect.~4.1, that the
axioms imply the following commutation relations between the
coefficients of twisted vertex operators:
\begin{equation}    \label{COMMUTATOR}
[A_{(m)}^{M^\sigma},B_{(k)}^{M^\sigma}]= \sum_{n\geq 0} \left(
                        \begin{array}{c} m \\ n \end{array} \right)
                        (A\nn \cdot B)^{M^\sigma}_{(m+k-n)},
\end{equation}
where by definition
$$
\left( \begin{matrix} m \\
n \end{matrix} \right) = \frac{m(m-1)\ldots (m-n+1)}{n!}, \quad n \in
\Z_{>0}; \qquad \left( \begin{matrix} m \\
0 \end{matrix} \right) = 1.
$$

We also have the following analogue of Prop.~4.1 of \cite{FBZ}:

\begin{lemma}    \label{TA}
For any $A \in V$, $Y^{\Mt}(TA,\zf{N}) = \pa_z Y^{\Mt}(A,\zf{N})$.
\end{lemma}

\begin{proof}
We apply axiom (2) in the situation where $B = \vac$. Then
$$
Y^{\Mt}(Y(A,z-w) \vac, \wf{N})v = \sum_{n\geq 0} Y^{\Mt}(A_{(-n-1)}
\vac,\wf{N})v (z-w)^n.
$$
But $A_{(-2)} \vac = TA$, therefore $Y^{\Mt}(TA,\wf{N})v$ appears as
the coefficient in front of $(z-w)$ in this series. Hence it should
coincide with the coefficient in front of $(z-w)$ in the expansion of
$Y^{\Mt}(A,\zf{N}) v$ in a power series in $\wf{N}$ and $(z-w)$. But
the latter is equal to $\pa_w Y^{\Mt}(A,\wf{N})$.
\end{proof}

Applying formula \eqref{COMMUTATOR} in the case when $A=\omega$
and $m=1$ (so that $A_{(m)} = L_0$), we obtain that
$$
[L^{\Mt}_0,B_{(k)}^{\Mt}] = (L_0 \cdot B)_{(k)}^{\Mt} + (L_{-1} \cdot
B)_{(k+1)}^{\Mt}.
$$
But in a conformal vertex algebra $L_{-1} \cdot B = T B$ and
$(TB)_{(k+1)}^{\Mt} = (-k-1) B_{(k)}$ by Lemma \ref{TA}. Therefore if
$B$ is homogeneous of degree $\Delta$, then
\begin{equation}    \label{comm with Bk}
[L^{\Mt}_0,B_{(k)}^{\Mt}] = (\Delta-k-1) B_{(k)}^{\Mt}.
\end{equation}

Suppose that $\Mt$ is an admissible module. Then we define a linear
operator $S_\sigma$ on $\Mt$ as follows. It acts on the eigenvectors
of $L^{\Mt}_0$ with eigenvalue $\frac{m}{N}$ by multiplication by
$e^{\frac{2\pi i m}{N}}$. Hence we obtain an action of the cyclic
group of order $N$ generated by $\sigma$ on $\Mt$, $g \mapsto S_g$.
According to the axioms of twisted module and formula \eqref{comm
with Bk} we have the following identity:
\begin{equation}    \label{action on module}
S_g^{-1} Y^{\Mt}(g \cdot A,\zf{N}) S_g = Y^{\Mt}(A,\zf{N}).
\end{equation}

Finally, we remark that there is an analogue of the Reconstruction
Theorem for twisted modules. Namely, suppose that $V$ is generated by
vectors $a^\al \in V, \al \in S$, in the sense of the usual
Reconstruction Theorem (see Theorem 4.5 of \cite{KAC} or Theorem 3.6.1
of \cite{FBZ}). Then if $M^\sigma$ is a $\sigma$--twisted $V$--module,
the twisted vertex operators $Y^{M^\sigma}(A,z^{\frac{1}{N}})$ for all
$A \in V$ may be reconstructed from the series
$Y^{M^\sigma}(a^\al,z^{\frac{1}{N}}), \al \in S$. This follows from
H.Li's formula for $Y^{M^\sigma}(A_{(n)} B,z^{\frac{1}{N}})$ in terms
of $Y^{M^\sigma}(A,z^{\frac{1}{N}})$ and
$Y^{M^\sigma}(B,z^{\frac{1}{N}})$ \cite{LI}. But this formula is more
complicated than its untwisted analogue, so the resulting formula for
a general twisted vertex operator usually looks rather cumbersome (see
for example formula \eqref{complicated} below).

\section{Torsors and twists} \label{TT}

Let $\Mt$ be an admissible conformal $\sigma_{V}$--twisted $V$--module
where $\on{ord}(\sigma_{V}) = N$.  In this section we define a group
$\AutOr{N}$ which naturally acts on $\Mt$, as well as natural torsors
for $\AutOr{N}$. This will allow us to twist $\Mt$ by a certain torsor
of formal coordinates.

\subsection{The group $\AutOr{N}$}

Let $\Aut \fps{\zf{N}}$ denote the group of continuous algebra
automorphisms of $\fps{\zf{N}}$. Since $\C[[z^{\frac{1}{N}}]]$ is
topologically generated by $z^{\frac{1}{N}}$, an automorphism $\rho$
of $\C[[z^{\frac{1}{N}}]]$ is completely determined by the image of
$z^{\frac{1}{N}}$, which is a series of the form
\begin{equation}
\rho(\zf{N}) = \sum_{n \in \frac{1}{N} \mathbb{Z}, n > 0}
c_{n} z^{n}, \label{rhodef1}
\end{equation}
where $c_{\frac{1}{N}} \ne 0$. Hence we identify $\Aut \fps{\zf{N}}$
with the space of power series in $\zf{N}$ having non-zero linear
term.  For more on the structure of the group $\Aut \fps{\zf{N}}$, see
Section 5.1 of \cite{FBZ}.  Recall that we denote $\C[[z]]$ by $\OO$.

\begin{definition} $\AutOr{N}$ is the subgroup of $\Aut \fps{\zf{N}}$
preserving the subalgebra $\fps{z} \subset \C[[z^{\frac{1}{N}}]]$.
\end{definition}

\noindent Thus, $\AutOr{N}$ consists of power series of the form
\begin{equation}
\rho(\zf{N}) = \sum_{n \in \frac{1}{N} + \mathbb{Z}, n > 0}
c_{n} z^{n}, \qquad c_{\frac{1}{N}} \ne 0. \label{rhodef}
\end{equation}

There is a homomorphism $\mu: \AutOr{N} \rightarrow \AutO$ which takes
$\rho \in \fps{\zf{N}}$ to the automorphism of $\fps{z}$ that it
induces. At the level of power series, this is just the map $\mu:
\rho(z) \mapsto \rho(z)^{N}$. The kernel consists of the automorphisms
of the form $\zf{N} \mapsto \epsilon \zf{N}$, where $\epsilon$ is an
$N$th root of unity, so we have the following exact
sequence:
\[
	1 \rightarrow \mathbb{Z}/ N \mathbb{Z} \rightarrow
	\AutOr{N} \rightarrow \AutO \rightarrow 1 \, .
\]
Moreover, $\AutOr{N}$ is a central extension of $\AutO$ by the cyclic
group $\mathbb{Z}/ N \mathbb{Z}$.

\noindent The Lie algebra of $\Aut \fps{\zf{N}}$ is $$\Der^{(0)}
\fps{\zf{N}} = \zf{N} \fps{\zf{N}} \partial_{\zf{N}},$$ and the Lie
algebra of $\AutOr{N}$ is its Lie subalgebra $\Der^{(0)}_{N} \OO =
\zf{N} \C[[z]] \partial_{\zf{N}}$. The homomorphism $\mu$ induces an
isomorphism of the corresponding Lie algebras sending
\[
	z^{k+\frac{1}{N}} \partial_{\zf{N}} \mapsto N z^{k}
	\partial_{z}, \qquad k \in \mathbb{Z}, k \geq 0.
\]

\subsection{The $\AutOr{N}$--torsor of special coordinates}
\label{torsor}

Let $(\mathcal{D},\sigma_{\mathcal{D}})$ be a pair consisting of a
formal disc $\mathcal{D} = \Spec R$, where $R \cong \fps{z}$ and an
automorphism $\sigma_{\mathcal{D}}$ of $\mathcal{D}$ (equivalently, of
$R$) of order $N$. We denote by $\bD$ the quotient of $\DD$ by
$\langle \sigma_{\DD} \rangle$, i.e., the disc $\Spec
R^{\sigma_{\DD}}$, where $R^{\sigma_{\DD}}$ is the subalgebra of
$\sigma_{\DD}$--invariant elements.

A formal coordinate $t$ is called a {\em special coordinate} with
respect to $\sigma$ if $\sigma(t) = \epsilon t$, where $\epsilon$ is
an $N$th root of unity, or equivalently, if $t^N$ is a formal
coordinate on $\bD$. We denote by $\AutD{}$ the set of all formal
coordinates on $\DD$ and by $\AutDN{N}{}$ the subset of $\AutD{}$
consisting of special formal coordinates. The set $\AutDN{N}{}$
carries a simply transitive right action of the group $\AutOr{N}$
given by $t \mapsto \rho(t)$, where $\rho$ is the power series given
in \eqref{rhodef}, i.e., $\AutDN{N}{}$ is an $\AutOr{N}$--torsor.

\subsection{Twisting modules by $\AutDN{N}{}$ }
\label{twistingmodules}

Let $\Mt$ be an admissible $\sigma_V$--twisted module over a conformal
vertex algebra $V$. Define a representation $r^{\Mt}$ of the Lie
algebra $\Der^{(0)}_N \OO$ on $\Mt$ by the formula
\[
	z^{k+\frac{1}{N}} \partial_{\zf{N}} \rightarrow - N \cdot
	L^{\Mt}_{k}.
\]
It follows from the definition of a twisted module that the operators
$L^{\Mt}_{k}, k > 0$, act locally nilpotently on $\Mt$ and that the
eigenvalues of $L^{\Mt}_{0}$ lie in $\frac{1}{N} \mathbb{Z}$, so that
the operator $N \cdot L^{\Mt}_{0}$ has integer eigenvalues. This
implies that the Lie algebra representation $r^{\Mt}$ may be
exponentiated to a representation $R^{\Mt}$ of the group $\AutOr{N}$.

In particular, the subgroup $\Z/N\Z$ of $\AutOr{N}$ acts on $\Mt$ by
the formula $i \mapsto S_\sigma^i$, where $S_\sigma$ is the operator
defined in \secref{algebras and modules}.

We now twist the module $\Mt$ by the action of $\AutOr{N}$ and define
the vector space
\begin{equation}
	\mathcal{\Mt}(\mathcal{D}) \overset{\on{def}}{=} \AutDN{N}{}
	\underset{\AutOr{N}}\times \Mt. \label{GeoTwistModule}
\end{equation}
Thus, vectors in $\mathcal{\Mt}({\mathcal{D}})$ are pairs $(t,v)$, up
to the equivalence relation $$(\rho(t),v) \sim (t,R^{\Mt}(v)), \qquad
t \in \AutDN{N}{}, v \in \Mt.$$ When $\fd = \fd{x}$, the formal
neighborhood of a point $x$ on an algebraic curve $X$, we will use the
notation $\mc{\Mt}_{x}$.

\section{Twisted vertex operators as sections}

Our goal is to give a coordinate-independent description of the
operation $Y^{\Mt}$. In order to do this we need to find how the
operation $Y^{\Mt}$ transforms under changes of special
coordinates. This is the subject of this section.

\subsection{The transformation formula for twisted vertex operators}
\label{TransFormSection}

Let $\OO = \C[[t]]$.  Denote by $R^{V}$ the representation of the
group $\Aut \OO$ on $V$ obtained by exponentiating the representation
$r^V$ of the Lie algebra $\Der^{(0)} \OO$ sending $t^{n+1} \pa_t$ to
$-L_n, n \geq 0$ (see Section 5.2 of \cite{FBZ}). Recall that for any
$\rho(t^{\frac{1}{N}}) \in \AutOr{N}$, we have
$\rho(t^{\frac{1}{N}})^N \in \Aut \OO$. For any $\tau(t) \in \AutO$ we
denote by $\tau_z$ the element of $\Aut (\C[[z]] \widehat{\otimes}
\OO)$ obtained by expanding $\tau(z+t) - \tau(z)$ in powers of $t$
(see Section 5.4.5 of \cite{FBZ}). Then we have the following analogue
of Lemma 5.4.6 from \cite{FBZ} (that lemma is originally due to
Y.-Z. Huang \cite{H}).

\begin{lemma} For any $A \in V$, $\rho \in \AutOr{N}$
\begin{equation} \label{TransformationFormula}
	R^{\Mt}(\rho) Y^{\Mt}(R^{V}((\rho^N)_{z})^{-1} A,
	\rho(\zf{N})) R^{\Mt}(\rho)^{-1} = Y^{\Mt} (A,\zf{N}).
\end{equation}
\end{lemma}

\begin{proof}
The exponential map $\Der^{(0)}_N \to \AutOr{N}$ is surjective, so it
suffices to consider the infinitesimal version of
\eqref{TransformationFormula}. Write
\[
	\rho = \exp{(\epsilon v(\zf{N}) \partial_{\zf{N}})} \cdot
	\zf{N},
\]
where
\[
	v(\zf{N})= - \sum_{k \in \mathbb{Z}, k \geq 0 } v_{k}
	z^{k+\frac{1}{N}}.
\]
We have
\[
	\rho^N = \exp{( \epsilon u(z) \partial_{z}) } \cdot z,
\]
where
\[ 
	u(z) = - N \sum_{k \in \mathbb{Z}, k \geq 0 } v_{k} z^{k+1}.
\]
To check that formula \eqref{TransformationFormula} holds, it suffices
to check that the $\epsilon$--linear term in it vanishes. Denote
$r^{\Mt}(v(\zf{N}) \partial_{\zf{N}})$ by $r^{\Mt}_{v}$ and the
$\ep$--linear term in $R^V((\rho^N)_z)$ by $r^{V}_{u,z}$.
The $\ep$--linear term in \eqref{TransformationFormula}
reads
\begin{equation*}
	(Id + \epsilon r^{\Mt}_{v}) Y^{\Mt}( (Id - \epsilon
	r^{V}_{u,z}) A, \zf{N} + \epsilon v(\zf{N}) ) (Id - \epsilon
	r^{\Mt}_{v}) - Y^{\Mt}(A,\zf{N})
\end{equation*}
\begin{equation}    \label{lin term}
= \epsilon [r^{\Mt}_{v},
	Y^{\Mt}(a,\zf{N})] - \epsilon Y^{\Mt}(r^{V}_{u,z} \cdot A,
	\zf{N}) + \epsilon v(\zf{N}) \partial_{\zf{N}}
	Y^{\Mt}(A,\zf{N}).
\end{equation}
We find that
\[
	r^{V}_{u,z} \cdot A = - \sum_{m \geq 0} \frac{1}{(m+1) !}
	(\partial^{m+1}_{z} u(z))L_{m} \cdot A
\]
and
\[
	r^{\Mt}_{v} = N \sum_{m \in \mathbb{Z}, m \geq 0} v_{m}
	L^{\Mt}_{m},
\]
so that vanishing of \eqref{lin term} is equivalent to the identity
\[
	[r^{\Mt}_{v}, Y^{\Mt}(A,\zf{N})] = - \sum_{m \geq -1}
	\frac{1}{(m+1)!} (\partial^{m+1}_{z} u(z)) Y^{\Mt}(L_{m} \cdot
	A, \zf{N}).
\]
Since
\[
	v(\zf{N}) \partial_{\zf{N}} = u(z) \partial_{z},
\]
this identity follows from the OPE between a twisted vertex operator
and the Virasoro field in the same way as in Section 5.2.3 of
\cite{FBZ}.
\end{proof}

\subsection{Example: primary fields}

Recall that a vector $A \in V$ is called a primary vector of conformal
dimension $\Delta$ if it satisfies
$$
L_n A = 0, \quad n>0; \qquad L_0 A = \Delta A.
$$
As shown in Lemma 5.3.4 of \cite{FBZ}, the corresponding vertex
operator $Y(A,z)$ transforms under coordinate changes as an
endomorphism-valued $\Delta$--differential on the punctured disc. Now
formula \eqref{TransformationFormula} implies an analogous
transformation formula for the corresponding twisted vertex operator
$Y^{\Mt}(A,\zf{N})$.

\begin{corollary}
Let $A \in V$ be a primary vector of conformal dimension
$\Delta$, and $\rho \in \AutOr{N}$. Then
\begin{equation} \label{TwistedPrimary}
R^{\Mt}(\rho) Y^{\Mt}(A, \rho(\zf{N}) ) R^{\Mt}(\rho)^{-1} \left(
\partial_{z} (\rho^{N}(\zf{N})) \right)^{\Delta} =
Y^{\Mt}(A,\zf{N}).
\end{equation}
%where by $\partial_{z}(\rho^{N}(\zf{N}) )$ we mean the series
%$z^{\frac{-(N-1)}{N}} \rho^{N-1}(\zf{N})$.
\end{corollary}

\section{Coordinate-independent interpretation of twisted vertex
  operators} \label{TwistSec1}

\subsection{Recollections from \cite{FBZ}}    \label{reminder}

Let $X$ be a smooth curve and $\Au_X$ the principal $\AutO$--bundle of
formal coordinates on $X$. The fiber of $\Au_X$ at $x \in X$ is the
$\AutO$--torsor $\Au_x$ of formal coordinates at $x$ (see Section 5.4
of \cite{FBZ} for details). Given a conformal vertex algebra $V$, set
\begin{equation}    \label{def of V}
\V = \V_X = \Au_X \underset{\AutO}\times V.
\end{equation}
This is a vector bundle whose fiber at $x \in X$ is the
$\Au_x$--twist of $V$,
$$
\V_x = \Au_x \underset {\AutO}\times V.
$$
i.e., the set of pairs $(z,A)$ where $z$ is a formal coordinate at $x$
and $A \in V$, modulo the equivalence condition $(\rho(z),A) \sim
(z,R^V(\rho) \cdot A)$ for $\rho \in \AutO$.

As explained in Chapter 5 of \cite{FBZ}, for any $x \in X$ the vertex
operation $Y$ gives rise to a canonical section $\mc{Y}_x$ of the dual
bundle $\V^*$ on the punctured disc ${\mc D}_x^\times$ with values in
$\on{End} \V_x$. Equivalently, we have a canonical linear map
$$
\mc{Y}^{\vee}_x: \Gamma({\mc D}_x^\times,\V \otimes \Omega_X) \to
\on{End} \V_x, \qquad s \mapsto \on{Res}_x \langle \Y_x,s \rangle.
$$
If we choose a formal coordinate $z$ at $x$ and use it to trivialize
$\V|_{{\mc D}_x}$, then $\mc{Y}^{\vee}_x$ is given by the formula
$$
\mc{Y}^{\vee}_x(A \otimes z^n dz) = A_{(n)}.
$$
Furthermore, the map $\mc{Y}^{\vee}_x$ factors through the quotient
$$
U(\V_x) = \Gamma({\mc D}_x^\times,\V \otimes \Omega_X)/\on{Im} \nabla.
$$
The latter is a Lie algebra and the resulting map $U(\V_x) \to
\on{End} \V_x$ is a Lie algebra homomorphism (see Section 8.2 of
\cite{FBZ}).

More generally, let $M$ be an admissible $V$--module. We attach to it
a vector bundle ${\mc M}$ on $X$ in the same way as above. The module
operation $Y^M$ then gives rise to a canonical section $\mc{Y}^M_x$ of
$\V^*|_{{\mc D}_x^\times}$ with values in $\on{End} {\mc
M}_x$. Equivalently, we have a canonical linear map
$$
\mc{Y}^{M,\vee}_x: \Gamma({\mc D}_x^\times,\V \otimes \Omega_X) \to
\on{End} {\mc M}_x,
$$
which factors through $U(\V_x)$ (see Section 6.3.6 of \cite{FBZ}).

In this section we obtain analogous results for twisted modules over
vertex algebras.

\subsection{The vector bundle $\V^{H}_{X}$ } \label{VXH}

Let $C$ be a smooth projective curve, and $H \subset Aut(C)$ a finite
group of automorphisms of $C$. Suppose furthermore that $V$ is a
conformal vertex algebra, and that $H$ acts on $V$ by conformal
automorphisms. The vector bundle $\V_{C}$ carries an $H$--equivariant
structure lifting the action of $H$ on $C$. It is given by
\begin{equation} \label{EqStruct}
	h \cdot (p,(A,z)) \overset{\on{def}}{=} (h(p),(h(A), z \circ
	h^{-1}))
\end{equation}
where $z \circ h^{-1}$ is the coordinate induced at $h(p)$ from
$z$. Let $X = C / H$ be the quotient curve, and $\nu: C \to X$ the
quotient map, ramified at the points where $H$ has non-trivial
stabilizers. Denote by $\ol{C} \subset C$ (resp. $\ol{X} \subset X$)
the complement of the ramification points (resp. branch points) of
$\nu$, and by $\ol{\nu}: \ol{C} \to \ol{X}$ the restriction of
$\nu$. Thus, $\ol{C}$ is an $H$--principal bundle over $\ol{X}$. The
action of $H$ on $\ol{C}$ is free, and $\V_{\ol{C}}$ descends to a
vector bundle $\V^{H}_{\ol{X}}$ on $\ol{X}$. More explicitly,
\begin{equation}    \label{wtV}
	\V^{H}_{\ol{X}} = \Au_{\ol{C}} \underset{\AutO \times
	H}{\times} V
\end{equation}
Here, $H$ acts on $\Au_{\ol{C}}$ by $h(p,z)=(h(p),z \circ h^{-1})$,
and this action commutes with the action of $\AutO$. The actions of
$H$ and $\AutO$ on $V$ commute because $H$ is a conformal automorphism
of $V$, and thus commutes with the Virasoro action.

%We will use the notation of Section \ref{Curves}. First, we observe
%that since the $\sigma_V(\omega) = \omega$ by our assumption, the
%actions of $\AutO$ and $\sigma_V$ on $V$ commute. Therefore the vector
%bundle $\V$ carries an automorphism, also denoted by $\sigma_V$. Using
%this automorphism, we define a $\sigma_C$--equivariant structure on
%the pull-back $\nu^*(\V)$ of $\V$ to $C$, i.e., an automorphism
%$\wt{\sigma}_C$ of the total space of $\nu^*(\V)$,
%\begin{equation}
%	\widetilde{\sigma}_{C}: (( z, A),y) \mapsto
%	((z,\sigma_{V}(A)), \sigma_{C}(y)), \label{eqstructure}
%\end{equation}
%extending the action of $\sigma_{C}$ on the base.

The vector bundle $\V^{H}_{\ol{X}}$ possesses a flat connection
$\nabla^{H}$. If $z$ is a local coordinate $x \in \ol{X}$,
$\nabla^{H}$ is given by the expression $d + L^{V}_{-1} \otimes dz$.

%so $\nu^*(\V)$ acquires the flat connection $\wt{\nabla} =
%\nu^{*}(\nabla)$, which commutes with $\wt{\sigma}_{C}$.

\subsection{Modules along $H$--orbits} \label{OrbitModule}

Let $x \in X$. Then every point $p \in \nu^{-1}(x)$ has a cyclic
stabilizer of order $N$, which we denote $H_{p}$. Each $H_{p}$ has a
canonical generator $h_{p}$, which corresponds to the monodromy of a
small loop around $x$. For a generic point $p$, $H_p = \{ e \}$ and we
set $h_p = e$. Suppose that we are given the following data:

\begin{enumerate}
	\item A collection of admissible $V$--modules $\{ M^{h_{p}}_{p}
	\}_{p \in \nu^{-1}(x)}$, one for each point in the fiber, such
	that $M^{h_{p}}_{p}$ is $h_{p}$--twisted.

        \item A collection
	of maps $S_{g,p,g(p)}: M^{h_{p}}_{p} \mapsto
	M^{h_{g(p)}}_{g(p)}$, $g \in H$, $p \in \nu^{-1}(x)$,
	commuting with the action of $\AutOr{N}$ and satisfying
	$$S_{gk,p,gk(p)} = S_{g,k(p),gk(p)} \circ S_{k,p,k(p)},$$
	$$S^{-1}_{g,p,g(p)} = S_{g^{-1},g(p),p},$$ and
	$$S_{g,p,g(p)}^{-1} Y^{M^{h_{g(p)}}_{g(p)}} (g \cdot A,z^{})
	S_{g,p,g(p)} = Y^{M^{h_{p}}_{p}}(A,z).$$
 
	\item If $g \in H_{p}$, then $S_{g,p,p} = S_{g}$, where
	$S_{g}$ is the operator defined in \secref{algebras and
	modules}.
\end{enumerate}

Given a collection $\{ M^{h_{p}}_{p} \}_{p \in \nu^{-1}(x)}$, we can
form the collection $\{ \mc{M}^{h_{p}}_{p}(\fd{p}) \}_{p \in
\nu^{-1}(x)}$, where $\mc{M}^{h_{p}}_{p}(\fd{p})$ is the
$\AutOr{N}$--twist of $M^{h_{p}}_{p}$ by the torsor of special
coordinates at $p$. 

Let 
\[
 \overline{\mc{M}_{x}} = \bigoplus_{p \in \nu^{-1}(x)}
 \mc{M}^{h_{p}}_{p}(\fd{p})
\]

This is a representation of $H$, where
$H$ acts as follows. If $A \in M^{h_{p}}_{p}$, $\zf{N}_{p}$ is a
special coordinate at $p$, and $g \in H$, then
\[
g \cdot (A, \zf{N}_{p}) = (S_{g,p,g(p)} \cdot A, \zf{N}_{p} \circ g^{-1}) 
\]
Note that this action is well-defined since the $S$--operators commute
with the action of $\AutOr{N}$. Now, let $\mc{M}_{x} =
(\overline{\mc{M}_{x}})^{H}$, the space of $H$--invariants of
$\overline{\mc{M}_{x}}$. The composition of the inclusion $\mc{M}_{x}
\to \overline{\mc{M}_{x}}$ and the projection $\overline{\mc{M}_{x}}
\to \mc{M}^{h_{p}}_{p}(\fd{p})$ is an isomorphism for all $p \in
\nu^{-1}(x)$. For $v_p \in \mc{M}^{h_{p}}_{p}(\fd{p})$, denote by
$[v_p]$ the corresponding vector in $\mc{M}_{x}$.  Note that for each
$(A,\zf{N}_{p})_p \in \mc{M}^{h_{p}}_{p}(\fd{p})$, and $g \in H$,
$[(A, \zf{N}_{p})_p] = [(S_{g,p,g(p)} \cdot A, \zf{N}_{p} \circ
g^{-1})_{g(p)}]$ in $\mc{M}_{x}$.

\begin{definition} \label{OrbitModuleDef}
We call $\mc{M}_{x}$ a \emph{$V$--module along $\nu^{-1}(x)$}.
\end{definition}

Henceforth, we will suppress the square brackets for
elements of $\mc{M}_{x}$ and refer to $[(A,\zf{N}_{p})_p]$ simply as
$(A,\zf{N}_{p})$.

\subsection{Construction of Modules along $H$--orbits} 
\label{ModuleConstruction}

In this section we wish to give a construction of a module along
$\nu^{-1}(x)$ starting with a point $p \in \nu^{-1}(x)$ and an
$h_p$--twisted module $M^{h_{p}}_{p}$. Note that when $H_{p}$ is
trivial, this just an ordinary $V$--module $M$.

Thus, suppose we are given $p \in \nu^{-1}(x)$, and an
$h_{p}$--twisted module $M^{h_{p}}_{p}$. Observe that the monodromy
generator at the point $g(p)$ is $h_{g(p)} = gh_{p}g^{-1}$, i.e. the
monodromies are conjugate.

\begin{enumerate}
	\item For $g \in H$, define the module
	$M^{gh_{p}g^{-1}}_{g(p)}$ to be $M^{h_{p}}_{p}$ as a vector
	space, with the $V$--module structure given by the vertex
	operator
\begin{equation} \label{NewStructure}
	Y^{M^{gh_{p}g^{-1}}_{g(p)}}(A,z) =
	Y^{M^{h_{p}}_{p}}(g^{-1} \cdot A,z)
\end{equation}
	It is easily checked that this equips
	$M^{gh_{p}g^{-1}}_{g(p)}$ with the structure of a
	$gh_{p}g^{-1}$--twisted module. Furthermore, if $g \in H_{p}$,
	this construction results in an $h_{p}$--twisted
	module isomorphic to $M^{h_{p}}_{p}$.

	\item Recall that $M^{h_{q}}$ is canonically isomorphic to
	$M^{h_{p}}_{p}$ as a vector space by the previous item. Thus,
	if $q \in \nu^{-1}(x)$, and $g(q) \ne q$, define
	$S_{g,q,g(q)}$ to be the identity map.

	\item If $g \in H_{q}$, then $g$ is
	conjugate to an element $g' \in H_{p}$. Define $S_{g,q,q} =
	S_{g',p,p}$ also using the canonical identification.

\end{enumerate}

It is easy to check that this construction is well-defined, and
satisfies the requirements of definition \ref{OrbitModuleDef}.

\begin{remark} 
If $H_{p}$ is trivial, and $M = V$, then for any $g \in H$, the new
module structure \ref{NewStructure} is isomorphic to the old one, and
so the resulting module $\mc{M}_{x}$ along $\nu^{-1}(x)$ is isomorphic
to $\mc{V}^{H}_{x}$, the fiber of the sheaf $\V^{H}$ at $x$.
\end{remark}

\begin{remark}
If $H=H_{p}$, then $p$ is unique, and so any $h_{p}$--twisted module
results in a module along $\nu^{-1}(x)$.
\end{remark}

%\begin{example}  \label{Example1}
%Suppose that we are given an ordinary (untwisted) $V$--module $M$
%together with an action of $H$ such that $g^{-1} Y^{M}(g\cdot A,z) g =
%Y^{M}(A,z), \forall g \in H$ (in particular, this holds for $M=V$). We
%show that if $x \in \ol{X}$ (i.e. the fiber over $x$ contains $\vert H
%\vert$ points), then $M$ gives rise to a $V$--module along
%$\nu^{-1}(x)$ with the following data: for $p \in \nu^{-1}(x)$, define
%$M^{e}_{p} = M$, and let $S_{g,p,g(p)} = g$.  The requirements of
%Definition \ref{OrbitModuleDef} are easily checked.
%\end{example}

%\begin{example}  \label{Example2}
%Suppose that $H$ is cyclic, generated by $h_p$, and that $\nu^{-1}(x)$
%consists of a single point $p$. Let $h_{p} \in H_{p} = H$ be as
%above. We show that an admissible $h_{p}$--twisted module $M^{h_{p}}$
%gives rise to a module along $\nu^{-1}(x)$. Let $M^{h_{p}}_{p} =
%M^{h_{p}}$, and set $S_{h,p,p} = S_{h}$. Recall from
%\secref{twistingmodules} that the assignment $h \mapsto S_{h}$
%coincides with the action of the subgroup $\Z/N\Z$ of $\AutOr{N}$ on
%$M^{h_{p}}$. The element $h \in H$ corresponds to the element of
%$\AutOr{N}$ given by the formula $h(\zf{N}) = \zf{N} \circ h^{-1}$.
%\end{example}

\subsection{Twisted vertex operators as sections of $\V^{H,*}_{X}$}
 \label{TvopsSec}
%Let $p$ be a ramification point on $C$. Formula
%\eqref{TransformationFormula} tells us that under the action of
%$\AutOr{N}$, the twisted vertex operation $Y^{\Mt}$ transforms as a
%kind of multi-valued section of the vector bundle
%$\V^{*}|_{\mathcal{D}_{\nu(p)}^\times}$. The monodromy of such a
%section may be resolved by pulling it back to a branched cover.

%We use the notation of Section \ref{TransFormSection}. Let
%$z^{\frac{1}{N}}$ be a special coordinate at $p$. Then the
%corresponding coordinate $z$ at $\nu(p)$ allows us to trivialize the
%bundle $\V$ on $\mathcal{D}_{\nu(p)}$. Let $\iota_{z}(A)$ denote the
%section of $\V|_{\mathcal{D}_{\nu(p)}}$ corresponding to $A$ in this
%trivialization. Denote ${\mc M}^\sigma({\mathcal{D}}_p)$ by ${\mc
%M}^\sigma_p$. Then we have the following analogue of Theorem 5.4.4 of
%\cite{FBZ}.

We begin with the observation that a section of $\V^{H}_{\ol{X}}$ over
$U \subset \ol{X}$ is the same as an $H$--invariant section of
$\V_{\ol{C}}$ over $\nu^{-1}(U)$, and likewise for
$\V^{H,*}_{\ol{X}}$. Thus, defining a section of $\V^{H,*}_{\ol{X}}$
on $\fpd{x}$ is equivalent to defining an $H$--invariant section of
$\V^{*}_{\ol{C}}$ on $\coprod_{p \in \nu^{-1}(x)} \fpd{p}$.

Let $p \in \nu^{-1}(x)$, and let $\zf{N}_{p}$ be an $h_{p}$--special
formal coordinate at $p$. This coordinate gives us a trivialization
$\iota_{z_{p}}$ of $\ds \left. \V_{\ol{C}} \right|_{\fpd{p}}$. We will
denote by $\iota_{z_{p}}(A)$ the section of $\ds \left. \V_{\ol{C}}
\right|_{\fpd{p}}$ corresponding to $A \in V$ with respect to this
trivialization. The coordinate $\zf{N}_{p}$ also gives us an
identification of $\mc{M}^{h_{p}}_{p}(\fd{p})$ with $M^{h_{p}}_{p}$.

%In what follows by an $\on{End}(\mc{M}_{x})$--valued section of
%$\V^{*}_{\ol{C}}$ on $\coprod_{p \in \nu^{-1}(x)} \fpd{p}$ we will
%understand a collection of
%$\on{End}(\mc{M}^{h_{p}}_{p}(\fd{p}))$--valued sections of
%$\V^{*}_{\ol{C}}$ on $\fpd{p}$ for all $p \in \nu^{-1}(x)$.

\begin{theorem} 	\label{Interpretation1}
	Let $x \in X$, and $\mc{M}_{x}$ a $V$--module along
	$\nu^{-1}(x)$. For each $p \in \nu^{-1}(x)$, choose an
	$h_{p}$--special coordinate $\zf{N}_{p}$ at $p$. Define an
	$\on{End}(\mc{M}_{x})$--valued section
	$\mc{Y}^{\mc{M}_{x}}$ of $\V^{*}_{\ol{C}}$ on
	$\coprod_{p \in \nu^{-1}(x)} \fpd{p}$ by the formula
\begin{equation}
	\langle (\zf{N}_{p},\phi),
	\mathcal{Y}^{\mc{M}_{x}}(\iota_{z_{p}}(A)) \cdot
	(\zf{N}_{p},v) \rangle = \langle \phi,
	Y^{M^{h_{p}}_{p}}(A,\zf{N}_{p}) \cdot v \rangle.
\end{equation}
Then this section $\mc{Y}^{\mc{M}_{x}}$ is independent of the choice
of special coordinate $\zf{N}_{p}$ on each $\fpd{p}$. Furthermore, it
is $H$--invariant.
\end{theorem}

%\begin{theorem}	\label{Interpretation1}
%	Define an $\End \mathcal{M}^\sigma_p$--valued section
%$\mathcal{Y}^{\Mt}_p$ of $\nu^{*}(\V)^{*}|_{\mathcal{D}^\times_p}$ by
%the formula
%\begin{equation}
%	\langle (\zf{N}_{p},\phi),
%	\mathcal{Y}^{\Mt}_p(\nu^{*}(\iota_{z}(A)))
%	\cdot (\zf{N}_{p},v) \rangle = \langle \phi, Y^{\Mt}(A,\zf{N}_{p})
%	\cdot v \rangle.
%\end{equation}
%Then this section $\mathcal{Y}^{\Mt}_p$ is independent of the choice
%of the coordinate $z^{\frac{1}{N}}$. Furthermore, this section is
%invariant with respect to the automorphism $\wt{\sigma}_C$ given by
%\eqref{eqstructure}.
%\end{theorem}

\begin{proof}

We begin by checking coordinate-independence. Choose a $p \in \nu^{-1}(x)$.
Let $w^{\frac{1}{N}}$ be another special coordinate at $p$. Then there
exists a unique $\rho \in \AutOr{N}$ such that $w^{\frac{1}{N}}=
\rho(\zf{N}_{p})$  (thus, $w = \rho^N(\zf{N}_{p})$). 
%Note also that $w$ and $z$ are formal coordinates at
%$x$ and $w = \rho^N(z)$. 
We attach to $w^{\frac{1}{N}}$ a section
$\widetilde{\mathcal{Y}}^{\Mt}_p$ of $\V_{\ol{C}}^{*}|_{\fpd{p}}$ with
values in $\on{End}(\mc{M}_{x})$ by the formula
\[
	\langle (w^{\frac{1}{N}},\widetilde{\phi}),
	\widetilde{\mc{Y}}^{\mc{M}_{x}}(\iota_{w}(\widetilde{A}))
	\cdot (w^{\frac{1}{N}},\widetilde{v}) \rangle = \langle
	\widetilde{\phi}, Y^{M^{h_{p}}}(\widetilde{A},w^{\frac{1}{N}}) \cdot
	\widetilde{v} \rangle.
\]
We must show that
$\widetilde{\mathcal{Y}}^{\mc{M}_{x}} = \mathcal{Y}^{\mc{M}_{x}}$. We have
\[
	(\zf{N}_{p},\phi) = (\wf{N},\phi \cdot R^{M^{h_{p}}}(\rho))
\] 
\[
	(\zf{N}_{p},v) = (\wf{N}, R^{M^{h_{p}}}(\rho)^{-1} \cdot v)
\]
As explained in \cite{FBZ}, the section $\iota_{z_{p}}(A)$ of
$\V_{\ol{C}}$ appears in the coordinate $w = \rho^{N}(\zf{N}_{p})$ as
$\iota_{w}(R^{V}((\rho^N)_{z})^{-1} \cdot A)$. It follows that
\[
	\iota_{z_{p}}(A) =
	\iota_{w} ( R^{V} ((\rho^N_{z})^{-1}) \cdot A ).
\]
Therefore
\begin{multline*}
	\langle (\zf{N}_{p},\phi),
	\widetilde{\mathcal{Y}}^{\mc{M}_{x}}(\iota_{z_{p}}(A)
	\cdot (\zf{N}_{p},v) \rangle \\ = \langle \phi, R^{M^{h_{p}}}(\rho)
	Y^{\Mt} (R^{V}((\rho^N)_{z})^{-1} \cdot A, \rho(\zf{N}_{p}))
	R^{M^{h_{p}}}(\rho)^{-1} v \rangle.
\end{multline*}
By \eqref{TransformationFormula}, $\widetilde{\mathcal{Y}}^{\mc{M}_{x}} =
\mathcal{Y}^{\mc{M}_{x}}$, and we obtain that our section is
coordinate-independent. \\

\noindent We now proceed to show that $\mc{Y}^{\mc{M}_{x}}$ is $H$--invariant. 
This amounts to checking, for $g \in H$ 
\[	
	\langle (\zf{N}_{p},\phi),
	\mc{Y}^{\mc{M}_{x}}(g \cdot \iota_{z_{p}}(A))
	\cdot (\zf{N}_{p},v) \rangle = g \cdot 
	\langle (\zf{N}_{p},\phi),
	\mc{Y}^{\mc{M}_{x}}(\iota_{z_{p}}(A))
	\cdot (\zf{N}_{p},v) \rangle 
\]
where the action on the right is by pullback of functions. The
right-hand side is $ \langle \phi, Y^{M^{h_{p}}_{p}}(A,\zf{N}_{p} \circ
g^{-1}) \cdot v \rangle $. On the left, we have
\begin{align*}
g \cdot \iota_{z_{p}}(A) &= \iota_{z_{p} \circ g^{-1}}(g \cdot A) \\
(\zf{N}_{p},\phi) & \cong (\zf{N}_{p} \circ g^{-1}, \phi \circ
S^{-1}_{g,p,g(p)}) \in \mc{M}^{*}_{x} \\ (\zf{N}_{p},v) & \cong
(\zf{N}_{p} \circ g^{-1}, S_{g,p,g(p)} \cdot v) \in \mc{M}_{x}
\end{align*}
Thus we get
\begin{align*}
	\langle (\zf{N}_{p},\phi),
	\mathcal{Y}^{\mc{M}_{x}} & (g \cdot
	\iota_{z_{p}}(A) \cdot (\zf{N}_{p},v) \rangle \\ &= \langle
	(\zf{N}_{p} \circ g^{-1}, \phi \circ S^{-1}_{g,p,g(p)}),
	\mc{Y}^{\mc{M}_{x}}( \iota_{z_{p} \circ g^{-1}}(g \cdot
	A)),(\zf{N}_{p} \circ g^{-1}, S_{g,p,g(p)} \cdot v) \rangle \\
	&= \langle \phi, S^{-1}_{g,p,g(p)} Y^{M^{h_{g(p)}}_{g(p)}}(g \cdot
	A,\zf{N}_{p} \circ g^{-1}) S_{g,p,g(p)} \cdot v \rangle \\ & =
	\langle \phi, Y^{M^{h_{p}}_{p}}(A,\zf{N}_{p} \circ g^{-1}) \cdot v
	\rangle \\
	&=  g \cdot 
	\langle (\zf{N}_{p},\phi),
	\mathcal{Y}^{\mc{M}_{x}}(\iota_{z_{p}}(A))
	\cdot (\zf{N}_{p},v) \rangle 
%g \cdot \langle \phi, Y^{M_{h_{p}}}(A,\zf{N}_{p}) \cdot v
%	\rangle 
\end{align*}
\end{proof}

\begin{remark}
  In view of the comments at the beginning of section \ref{TvopsSec},
  $\mc{Y}^{\mc{M}_{x}}$ is the pullback under $\nu$ of a unique section
  of $\vbund$ on $\fpd{x}$. We will abuse notation by
  denoting the latter by $\mc{Y}^{\mc{M}_{x}}$ as well.\qed
\end{remark}

%The statement of the theorem may be interpreted as follows. Let
%$\ol{C}$ (resp., $\ol{X}$) be the complement of the ramification
%points in $C$ (resp., $X$). Then $\ol{\nu}: \ol{C} \to \ol{X}$ is a
%principal $\Z/N\Z$--bundle. The group $\Z/N\Z$ acts on the bundle
%$\V_{\ol{X}}$ through the powers of the automorphism $\sigma_V$. Let
%\begin{equation}    \label{wtV}
%\wt{\V}_{\ol{X}} = \ol{C} \underset{\Z/N\Z}\times \V_{\ol{X}}
%\end{equation}
%be the twist of $\V_{\ol{X}}$ by the bundle $\Z/N\Z$. Then
%Theorem~\ref{Interpretation1} may be reformulated as

%\begin{corollary}    \label{torsor1}
%The operation $Y^{\Mt}$ gives rise to an $\End
%\mathcal{M}^\sigma_p$--valued section of $\wt{\V}_{\ol{X}}^*$ over
%$\fpd{\nu(p)}$.
%\end{corollary}

In the case of twisted primary fields, \eqref{TwistedPrimary} implies
the following analogue of Proposition 5.3.8 of \cite{FBZ}:

\begin{proposition} \label{PrimaryField}
	Let $x \in X$, and $\mc{M}_{x}$ a $V$--module along
	$\nu^{-1}(x)$. For each $p \in \nu^{-1}(x)$, choose an
	$h_{p}$--special coordinate $\zf{N}_{p}$ at $p$. Define an
	$\on{End}(\mc{M}_{x})$--valued $\Delta$--differential $\varpi$ 
	on $\coprod_{p \in \nu^{-1}(x)} \fpd{p}$ by the formula
\begin{align*}
	\langle (\zf{N}_{p}, \phi), \varpi \cdot (\zf{N}_{p},v)
	\rangle &= \langle \phi, Y^{M^{h_{p}}_{p}}(A, \zf{N}_{p})
	\cdot v \rangle (dz_{p})^{\Delta} \\ &= N^{\Delta} z^{\Delta
	\frac{(N-1)}{N}}_{p} \langle \phi,
	Y^{M^{h_{p}}_{p}}(A,\zf{N}_{p}) \cdot v \rangle
	(d\zf{N}_{p})^{\Delta}
\end{align*}
Then $\varpi$ is independent of the choice of $\zf{N}_{p}$'s. 
\end{proposition}

Recall from Section 5.4.9 of \cite{FBZ} that a primary vector $A \in
V$ determines a line subbundle $j_{A} : \Omega_C^{-\Delta}
\hookrightarrow \V_{C}$, and by dualizing a surjection $j^{*}_{A}:
\V^{*}_{C} \twoheadrightarrow \Omega_C^{\Delta}$. The section $\varpi$
of $\Omega^{\Delta}_{\ol{C}}|_{\coprod_{p \in \nu^{-1}(x)} \fpd{p}}$
appearing in Proposition \ref{PrimaryField} is just the image of the section
$\mc{Y}^{\mc{M}_{x}}$ under $j^*_A$.

\subsection{Dual version}    \label{dual}

Let $x \in X$, and $\mc{M}_{x}$ a $V$--module along $\nu^{-1}(x)$.  As
in the case of ordinary vertex operators (see Section 5.4.8 of
\cite{FBZ}), dualizing the construction we obtain a linear map

$$
\mathcal{Y}^{\mc{M}_{x}, \vee}: \Gamma(\fpd{x},\V^{H}_{\ol{X}} \otimes
\Omega_{\ol{X}}) \to \End {\mc{M}_{x}}.
$$
Given by 
\begin{equation} \label{TwistLieAlg}
s \rightarrow \Res_{x} \langle \mc{Y}^{\mc{M}_{x}}, s \rangle 
\end{equation}
Moreover, this map factors through the quotient
$$ \label{PointAlgebra}
U(\V^{H}_{x}) \stackrel{\on{def}}{=} \Gamma(\fpd{x},\V^{H}_{\ol{X}}
\otimes \Omega_{\ol{X}})/\on{Im}
\nabla^{H},
$$
which has a natural Lie algebra structure. The corresponding map
$U(\V^{H}_{x}) \to \End {\mc{M}_{x}}$ is a homomorphism of Lie
algebras. Note that $x$ does not have to lie in $\ol{X}$, but can be
\emph{any} point of $X$.

\subsection{A sheaf of Lie algebras}    \label{sheaf of la}

Following Section 8.2.5 of \cite{FBZ}, let us consider the following
complex of sheaves (in Zariski topology) on $\ol{X}$:
$$
0 \rightarrow \V^{H}_{\ol{X}} \xrightarrow{\nabla}
\V^{H}_{\ol{X}} \otimes \Omega_{\ol{X}} \rightarrow 0
$$
where $\V^{H}_{\ol{X}} \otimes \Omega_{\ol{X}}$ is placed in cohomological
degree $0$ and $\V^{H}_{\ol{X}}$ is placed in cohomological degree
$-1$ (shifted de Rham complex). Let $h(\V^{H}_{\ol{X}})$ denote the sheaf of
the 0th cohomology, assigning to every Zariski open subset $\Sigma
\subset \ol{X}$ the vector space
$$
U_{\Sigma}(\V^{H}_{\ol{X}}) \overset{\on{def}}{=} \Gamma(\Sigma,\V^{H}_{\ol{X}}
\otimes \Omega_{\ol{X}})/\on{Im} \nabla^{H} 
$$
One can show as in Chapter 18 of \cite{FBZ} that this is a sheaf of
Lie algebras.
 
According to formula \eqref{PointAlgebra}, for any $x \in \Sigma'$,
where $\Sigma' \subset X$ is such that $\Sigma' \cap \ol{X} = \Sigma$,
restriction induces a Lie algebra homomorphism
$U_{\Sigma}(\V^{H}_{\ol{X}}) \to U(\V^{H}_{x})$. We denote the image
by $U_{\Sigma}(\V^{H}_{x})$.

\subsection{Interpretation in terms of chiral algebras}

A. Beilinson and V. Drinfeld have introduced in \cite{BD} the
notion of chiral algebra (see also \cite{Gai}). A chiral algebra on a
smooth curve $X$ is a right ${\mc D}$--module $\cA$ on $X$ together
with homomorphisms of ${\mc D}$--modules $\Omega \to \cA$ and
$$
j_* j^*(\cA \boxtimes {\cA}) \to \Delta_!({\cA }),
$$
where $\Delta: X \to X^2$ is the diagonal embedding, and $j: (X^2
\backslash \Delta) \to X^2$ is the complement of the diagonal. These
homomorphisms must satisfy certain axioms.

As shown in Chapter 18 of \cite{FBZ}, for any conformal vertex algebra
$V$ and any smooth curve $X$, the right ${\mc D}$--module $\V_X
\otimes \Omega_X$ is naturally a chiral algebra.

Recall that a module over a chiral algebra $\cA$ on $X$ is a right
${\mc D}$--module ${\mc R}$ on $X$ together with a homomorphism of
${\mc D}$--modules
$$
a: j_* j^*(\cA \boxtimes {\mc R}) \to \Delta_!({\mc R}).
$$
This homomorphism should satisfy the axioms of \cite{BD}.

Suppose that ${\mc R}$ is supported at a point $x \in X$ and denote
its fiber at $x$ by ${\mc R}_x$. Then ${\mc R} = i_{x!}({\mc R}_x)$,
where $i_x$ is the embedding $x \to X$. Applying the de Rham functor
along the second factor to the map $a$ we obtain a map
$$
a_x: j_{x*} j^*_x(\cA) \otimes {\mc R}_x \to {\mc R},
$$ where $j_x: (X \backslash x) \to X$. The chiral module axioms may
be reformulated in terms of this map. Note that it is not necessary
for $\cA$ to be defined at $x$ in order for this definition to make
sense. If $\cA$ is defined on $\X\backslash x$, we simply replace
$j_{x*} j^*_x(\cA)$ by $j_{x*}(\cA)$.

If $M$ is a module over a conformal vertex algebra $V$, we associate
to it the space ${\mc M}_x$ as in \secref{reminder}. Then the right
${\mc D}$--module $i_{x!}({\mc M}_x)$ is a module over the chiral
algebra $\V \otimes \Omega_X$ supported at $x$ and the corresponding
map $a^M_x$ is defined as follows. Choose a formal coordinate $z$ at
$x$ and use it to trivialize $\V|_{{\mc D}_x}$ and ${\mc M}_x$ and to
identify
$$
i_{x!}({\mc M}_x) = M((z))dz/M[[z]] dz.
$$
Then
$$
a^M_x(\imath_z(A) \otimes f(z)dz,B) = Y^M(A,z) B \otimes f(z)dz \quad
\on{mod} \quad M[[z]] dz.
$$
The independence of $a^M_x$ on $z$ is proved in the same way as the
independence of $\Y^M_x$. Note that applying to $a^M_x$ the de Rham
cohomology functor we obtain the map $\Y^{M,\vee}_x$ (see
\secref{reminder}).

Now let $\V^H_{\ol{X}}$ be the vector bundle on $\ol{X}$ with
connection defined by formula \eqref{wtV}. Then $\V^H_{\ol{X}} \otimes
\Omega_X$ is a right ${\mc D}$--module on $\ol{X}$. Suppose that $H =
\Z/N\Z = \langle \sigma \rangle$. Given a full ramification point $p
\in C$ and a twisted $V$--module $M^\sigma$, we define a map
$$
a^{\Mt}_p: j_{x*} (\V^H_{\ol{X}} \otimes \Omega_X)
\otimes {\mc M}^\sigma_p \to i_{p!}({\mc M}^\sigma_p),
$$
where $x=\nu(p) \in X \backslash \ol{X}$, as follows. Observe that
sections of $\V^H_{\ol{X}} \otimes \Omega_X|_{{\mc D}^\times_x}$
are the same as $\sigma$--invariant sections of $\nu^*(\V_X)
\otimes \Omega_C|_{{\mc D}^\times_p}$. Choose a special coordinate
$z^{\frac{1}{N}}$ at $p$ and use it to trivialize $\nu^*(\V)
\otimes \Omega_C|_{{\mc D}^\times_p}$ and ${\mc M}_p$ and to identify
$$
i_{p!}({\mc M}^\sigma_p) =
M((z^{\frac{1}{N}}))dz^{\frac{1}{N}}/M[[z^{\frac{1}{N}}]]
dz^{\frac{1}{N}}.
$$
Then
$$
a^{\Mt}_p(\imath_z(A) \otimes f(z^{\frac{1}{N}})dz^{\frac{1}{N}},B)
\overset{\on{def}}{=} Y^{\Mt}(A,z^{\frac{1}{N}}) B \otimes
f(z^{\frac{1}{N}})dz^{\frac{1}{N}} \quad \on{mod} \quad
M[[z^{\frac{1}{N}}]] dz^{\frac{1}{N}}.
$$
The independence of $a^{\Mt}_p$ on $z^{\frac{1}{N}}$ may be proved in
the same way as the independence of $\Y^{\Mt}_p$ was proved. Applying
to $a^{\Mt}_p$ the de Rham cohomology functor we obtain the map
$\Y^{\Mt,\vee}_p$ from \secref{dual}.

Now let $H$ be an arbitrary finite group acting (generically with
trivial stabilizers) on a smooth curve $C$. Then as before we have an
$H$--torsor $\ol{C}$ over $\ol{X} \subset X = C/H$. Let $\cA$ be a
chiral algebra on $\ol{X}$ equipped with an action of $H$ by
automorphisms. Then the $\ol{C}$--twist of $\cA$,
$$
\cA^{\ol{C}} = \ol{C} \underset{H}\times \cA,
$$
inherits the chiral algebra structure from $\cA$. So we can consider
$\cA^{\ol{C}}$--modules supported at arbitrary points $x \in X$. If
$\cA = \V_X \otimes \Omega_X$, where $V$ is a conformal vertex algebra
on which $H$ acts by automorphisms, then such modules may be
constructed from twisted $V$--modules. Namely, to each $V$--module
along $\nu^{-1}(x)$ (see Definition \ref{OrbitModuleDef}) we can
attach to it in the same way as above a $\cA^{\ol{C}}$--module
supported at $x$.

\section{Conformal blocks} \label{CBSec}

%Let $\{ p_{i} \}_{i=1,\dots,r}$ be the set of all $\sigma_{C}$--fixed
%points in $C$. Let $\{ q_{j} \}_{j=1,\ldots,s}$ be a collection of
%points such that $\# \Orb{q_j} = N$ and $\Orb{q_{j_{1}}} \cap
%\Orb{q_{j_{2}}} = \nul$ for $j_1 \neq j_2$. Let $\{ \Mt_{i}
%\}_{i=1,\ldots,r}$ be a collection of admissible $\sigma_{V}$--twisted
%$V$--modules and $\{ M_{j} \}_{j=1,\ldots,s}$ a collection of
%untwisted admissible $V$--modules.  We think of $\Mt_{i}$ as inserted
%at $p_{i}$, and $M_{j}$ as attached to the orbit $\Orb{q_{j}}$.  We
%will denote by $\Mtp{i}$ the twist of $\Mt_{i}$ by $\AutDN{N}{p_{i}}$,
%and by $\Mp{j}$ the twist of $M_{j}$ by $\AutD{\nu(q_{j})}$.

%We have seen that $Y^{\Mt_{i}}$ gives rise to an invariant $\End
%\Mtp{i}$--valued section $\mathcal{Y}^{M_i}_{p_i}$ of $\nV^{*}$ over
%${\fpd{p_{i}}}$, and $Y^{M_{j}}$ gives rise to an invariant $\End
%\Mp{j}$--valued section $\mathcal{Y}^{M_j}_{\Orb{q_{j}}}$ of
%$\nV^{*}$ over $\bigcup_k \fpd{\sigma_{C}^{k}(q_{j})}$. Let
%\[
%	\ol{C} = C \backslash (\{ p_{i} \} \cup \{ \Orb{q_{j}} \} ) 
%\]

We use the notation of section \ref{VXH}. Let $\{ x_{i} \}_{i=1 \cdots m}$ be
a collection of distinct points of $X$, which contains all of the
branch points of $\nu$. Let $\{\mc{M}_{x_{i}} \}_{i=1 \cdots m}$ be a
collection of $V$--modules, such that $\mc{M}_{x_{i}}$ is a module
along $\nu^{-1}(x_{i})$. Let
\[
	\mathcal{F} = \bigotimes^{m}_{i=1} \mc{M}_{x_{i}} 
\]
Given $\phi \in \mathcal{F}^{*}$, $A_{i} \in \mc{M}_{x_{i}}$,
\begin{equation}  \label{twists}
	\langle \phi, A_{1} \otimes \cdots \otimes
		\mathcal{Y}^{\mc{M}_{x_{i}}} \cdot A_{i} \otimes
	 \cdots \otimes A_{m} \rangle 
\end{equation}
is a section of $\vbund$ on $\fpd{x_{i}}$. 

We can now define the generalized space of conformal blocks, extending
Definition 9.1.1 of \cite{FBZ}:

\begin{definition}	\label{ConfB2}
The space of conformal blocks $$\Cvb{V}{m}$$ is by definition the
vector space of all linear functionals $\phi \in \mathcal{F}^{*}$ such
that for any $A_{i} \in \mc{M}_{x_{i}}$, the sections
\ref{twists} can be extended to the
\emph{same} section of $\vbund$, regular over 
$X \backslash \{x_{i} \}$. 
\end{definition}

\medskip

\noindent {\bf{Remark.}}  Observe that in this definition \emph{all}
branch points of $\nu$ are required to carry module insertions, so
that $X \backslash \{x_{i} \} \subset \ol{X}$. \qed

\medskip

We can pull back the sections \ref{twists} by $\nu$ to sections
of $\V^{*}_{\ol{C}}$ on $\coprod_{p \in \nu^{-1}(x_{i})} \fpd{p}$
Equivalently, Definition \ref{ConfB2} can be rephrased as follows:

\begin{definition}
The space of conformal blocks $$\Cvb{V}{m}$$ is by definition the
vector space of all linear functionals $\phi \in \mathcal{F}^{*}$ such
that for any $A_{i} \in \mc{M}_{x_{i}}$, the pullbacks of the sections
\ref{twists} can be extended to the
\emph{same} $H$--invariant section of $\V^{*}_{\ol{C}}$, regular over 
$C \backslash \{ \nu^{-1}(x_{i}) \}$. 
\end{definition}

%We can ask how the space of conformal blocks changes if we replace the
%points $\{ q_{j} \}$ with $\{ q'_{j} \}$, where $q'_{j} \in
%\Orb{q_{j}}$. Not surprisingly, the two spaces are isomorphic.

%\begin{lemma} If $q'_{j} \in \Orb{q_{j}}$ for all $j=1,\ldots,s$,
%then there is a canonical isomorphism
%\[
%	\Cvb{r}{s} \cong \Cvbtwo{r}{s}.
%\]
%\end{lemma}

\subsection{Alternative definition}

Composing the map \eqref{TwistLieAlg} with
the map
\[
U_{X \backslash \{ x_{i} \} } (\vbun) \rightarrow \bigoplus^{m}_{i=1}
U(\V^{H}_{x_{i}})
\]
we obtain an action of the Lie algebra $U_{X \backslash \{ x_{i}
\}}(\V^{H})$ on $\mc{F}$.  We will employ the Strong Residue Theorem
(see \cite{T}):

\begin{theorem} \label{ResThm}
Let $\mathcal{E}$ be a vector bundle on a smooth projective curve
$Z$. Let $t_{1},\ldots,t_{n} \in Z$ be a set of distinct
points. Then a section
\[
	\tau \in \bigoplus^{n}_{i=1} \Gamma(\fpd{t_{i}}, \mathcal{E}^{*})
\]
has the property that
\[
	\sum^{n}_{i=1} \Res_{t_{i}} \langle \mu, \tau \rangle =0,
		\qquad \forall \mu \in \Gamma(Z - \{ t_{1}, \cdots,
		t_{n}\}, \mathcal{E} \otimes \Omega_Z)
\]
if and only if $\tau$ can be extended to a regular section of
$\mathcal{E}$ over $Z - \{ t_{1}, \cdots, t_{n} \}$ (i.e., $\tau \in
\Gamma(Z - \{ t_{1}, \cdots , t_{n} \}, \mathcal{E}^{*})$) 
\end{theorem}

Applying Theorem \ref{ResThm} to Definition \ref{ConfB2}, we obtain
that $\phi \in \mc{F}^{*}$ is a conformal block if and only if it
vanishes on all elements of the form $s \cdot v,\; s \in U_{X
\backslash \{ x_{i} \}}(\vbun), v \in \mc{F}$. This leads to the
following equivalent definition, extending Definition 9.1.2 of
\cite{FBZ}:

\begin{definition}    \label{def coinv}
The space of \emph{coinvariants} is the vector space 
$$\Coinvb{V} = \mc{F} / U_{X \backslash \{x_{i} \}}(\V^{H}_{\ol{X}})
\cdot \mc{F}.
$$ The space of \emph{conformal blocks} is its dual: the vector space
of $U_{X \backslash \{x_{i} \} }(\V^{H}_{\ol{X}})$--invariant
functionals on $\mathcal{F}$ $$\Cvb{V}{m} = \on{Hom}_{U_{X \backslash
\{ x_{i} \} }(\V^{H}_{\ol{X}})}(\mc{F}, \mathbb{C}).$$
\end{definition}

\section{Example: Heisenberg vertex algebra} \label{HeisCB}

The definition of conformal blocks given in Section \ref{CBSec} is quite
abstract, and involves a priori all of the fields of the vertex
algebra $V$. When vertex algebras are generated by a finite number of
fields the definition of conformal block can be simplified to involve
only those generating fields.  In this section we illustrate this in
the case of the Heisenberg vertex algebra with an order $2$
automorphism.

\subsection{The vertex algebra $\pi$ and its $\mathbb{Z}/2
\mathbb{Z}$--twisted sector}

Let $\H$ (resp. $\Hs$) denote the Lie algebra with generators $\{
\widetilde{b}_{n}, \widetilde{K} \}_{n \in \mathbb{Z}}$ (resp. $\{
b_{n}, K \}_{n \in \frac{1}{2} + \mathbb{Z}}$ ), and commutation
relations
\[
	[ \widetilde{b}_{n}, \widetilde{b}_{m} ] = n \delta_{n,-m}
	\widetilde{K}
\] 
(resp. same but with $b_{n}$'s) where $\widetilde{K}$ (resp. $K$) is
central. Let $\H_{+}$ (resp. $\Hs_{+}$) be the subalgebras generated
by $\{ \widetilde{b}_{n} \}_{n \geq 0}$ (resp. $\{ b_{n} \}_{n >
0}$). For $\lambda \in \mathbb{C}$, let
$\widetilde{\mathbb{C}^{\lambda}}$ denote the $1$ dimensional
representation of $\H_{+} \oplus \mathbb{C} \cdot \wt{K}$ on which
$\widetilde{b}_{n} , n > 0$ acts by $0$, $\widetilde{b}_{0}$ acts by
$\lambda$, and $\widetilde{K}$ acts by the identity. Let
\[
	\pi^{\lambda} = \Ind^{\H}_{\H_{+} \oplus \mathbb{C} \cdot \wt{K}}
	\widetilde{\mathbb{C}^{\lambda}}.
\]
It is well-known that $\pi=\pi^{0}$ has the structure of a vertex
algebra (see for instance Chapter 2 of \cite{FBZ}), generated by the
field assignment
\[
	Y^{\pi}(\widetilde{b}_{-1} \vac, z ) = \widetilde{b}(z) =
	\sum_{n \in \mathbb{Z}} \widetilde{b}_{n} z^{-n-1}.
\]	
Let us take $\omega = \frac{1}{2} \widetilde{b}^{2}_{-1}$ to be the
conformal vector of $\pi$.  With this conformal structure, $\pi$ has a
conformal automorphism $\sigma$ of order $2$, induced from the
automorphism of $\H$ which acts by $\widetilde{b}_{n} \rightarrow -
\widetilde{b}_{n}$. All $\pi^{\lambda}$ have the structure of
conformal $\pi$--modules. If $\lambda = \sqrt{M}$, M even, then
$\pi^{\lambda}$ is admissible. We write
\[
	Y^{\pi^{\lambda}}(\widetilde{b}_{-1} \vac, z ) =
	\widetilde{b}^{\lambda}(z) = \sum_{n \in
	\mathbb{Z}} \widetilde{b}_{n} z^{-n-1}.
\]

Now let $\mathbb{C}$ denote the $1$--dimensional representation of
$\Hs_{+} \oplus \mathbb{C} \cdot K$ on which $K$ acts by the
identity. The $\Hs$--module
\[
	\pi^{\sigma} = \Ind_{\Hs_{+} \oplus \mathbb{C} \cdot K}^{\Hs}
	\mathbb{C}
\]
has the structure of an admissible conformal $\sigma$--twisted
$\pi$--module, generated (in the sense of \cite{LI}) by the field assignment
\[
	Y^{\pi^{\sigma}}(\widetilde{b}_{-1} \vac, z) = b(\zf{2}) =
	\sum_{n \in \frac{1}{2} + \mathbb{Z}} b_{n} z^{-n-1}.
\]
The twisted vertex operator assigned to an arbitrary vector $v \in
\pi$ is given as follows (see \cite{FLM,KP,Dong}). Let
\[
	W^{\pi^{\sigma}}(\widetilde{b}_{n_1} \ldots \wt{b}_{n_k} \vac, z) =
	\frac{1}{(-n_1-1)!} \ldots \frac{1}{(-n_k-1)!} \Wick
	\pa_z^{-n_1-1} b(z) \ldots \pa_z^{-n_k-1} b(z) \Wick
\]
and set
\[
	\Delta_{z} = \sum_{m,n \geq 0} c_{mn} \wt{b}_{m} \wt{b}_{n} z^{-m-n}
\] 
where the constants $c_{mn}$ are determined by the formula
\[
	\sum_{m, n \geq 0} c_{mn} x^{m} y^{n} = -
	\log \left( \frac{(1+x)^{1/2} + (1+y)^{1/2}}{2} \right)
\]
Then for any $v \in \pi$ we have
\begin{equation}    \label{complicated}
Y^{\pi^{\sigma}}(v, z) = W^{\pi^{\sigma}}( \exp{\Delta_{z}} \cdot v, z).
\end{equation}

\subsection{The Lie algebra ${\mc H}_{\on{out}}$} \label{Houtc}

Using the notation of earlier sections, we now restrict to the case
where $V = \pi$, $H \cong \langle \sigma_{C} \rangle$, where
$\sigma_{C}$ has order 2, and $\nu: C \rightarrow X$ has degree
$2$. This is for example the case when $C$ is hyperelliptic,
$\sigma_{C}$ is the hyperelliptic involution, and $X = \mathbb{C}
\mathbb{P}^{1}$. The vector bundles $\V_{\ol{C}}, \V^{H}_{\ol{X}}$
will be denoted $\Pi_{\ol{C}}, \Pi^{\sigma}_{\ol{X}}$ respectively.  

For $x \in X$ we can construct $\pi$--modules along $\nu^{-1}(x)$ by
applying the construction in section \ref{ModuleConstruction}. If
$\pi^{-1}(x)$ consists of one point, we obtain a module
$\pi^{\sigma}_{x}$ along $\nu^{-1}(x)$ starting with $\pi^{\sigma}$.
If $\pi^{-1}(x)$ consists of two points, we obtain a module
$\pi^{\lambda,p}_{x}$ along $\nu^{-1}(x)$ starting with a point $p$ in
the fiber and a $\pi^{\lambda}, \lambda \in \sqrt{2 \Z}$.

Let $\{ x_{i} \}_{i=1 \cdots m}$ be a collection of points of $X$
containing all of the branch points of $\nu$, and $\{ \pi_{x_{i}} \}$ a
collection of $\pi$--modules along $\nu^{-1}(x_{i})$, where $\pi_{x_{i}} \cong
\pi^{\sigma}_{x_{i}}$ or $\pi_{x_{i}} \cong \pi^{\lambda_{i},
p_{i}}_{x_{i}}$ depending on whether $\nu^{-1}(x_{i})$ consists of one
or two points. The vector $\widetilde{b}_{-1} \vac \in \pi$ is primary and has
conformal weight $1$.  Applying Proposition \ref{PrimaryField}, for each
$x_{i}$, we obtain an $\End(\pi_{x_{i}})$--valued $1$--form
$\varpi_{i}$ on $\coprod_{p \in \nu^{-1}(x_{i})} \fpd{p}$
%Theorem \ref{PrimaryField} implies that if
%$\zf{N}_{i}$ is a special coordinate centered at $p_{i} \in C$, then
%\[
%	\varpi_{i} = b(\zf{2}_{i}) d z = 2 \zf{2}_{i} b(\zf{2}_{i}) d
%	\zf{2}_{i}
%\]
%is a canonical $\End(\pi^{\sigma}_{i})$--valued one-form on $\fpd{p_{i}}$. 
%Similarly, let $z_{j}$ be a formal coordinate centered on
%$\nu(q_{j})$. Then $z_{j}$ induces a coordinate $w_{j} = z_{j} \circ
%\nu$ near the two points $q_{j}, \widetilde{q}_{j} =
%\sigma_{C}(q_{j})$. Define an $\End(\pi^{\lambda_{j}}_{j})$--valued
%one-form $\varpi_{j}$ on $\fpd{q_{j}} \bigcup \fpd{\widetilde{q}_{j}}$
%given by $\widetilde{b}^{\lambda_{j}}(w_{j}) dw_{j}$ on $\fpd{q_{j}}$,
%and $- \widetilde{b}^{\lambda_{j}}(w_{j}) dw_{j}$ on
%$\fpd{\widetilde{q}_{j}}$. The following observation is immediate: 
 %\begin{lemma}
%	The one-form $\varpi_{i}$ is the projection of the
%	$\on{End}(\pi^{\sigma}_{i})$--valued section
%	$\mathcal{Y}^{\pi^{\sigma}_{i}}$, and $\varpi_{j}$ is a
%	projection of the $\End(\pi^{\lambda_{j}}_{j})$--valued
%	section $\mathcal{Y}^{\pi^{\lambda_{j}}_{j}}_{q_{j}}$. 
%\end{lemma} 
%Let
%\[
%	\mathcal{F} = \bigotimes^{r}_{i=1} \pi^{\sigma}_{i} \otimes
%	\bigotimes^{s}_{j=1} \pi^{\lambda_{j}}_{j},
%\]
Let 
\[
   \mathcal{F} = \bigotimes^{m}_{i=1} \pi_{x_{i}} 
\]
and let $\Houtc$ be the abelian Lie algebra $\mathbb{C}[\Cnot ]$ of
  regular functions on $\Cnot$. If $A_{i} \in \pi_{x_{i}}$, then $f
  \in \Houtc$ acts on $\mc{F}$ by
%If $A_{i} \in \pi^{\sigma}_{i},
%  B_{j} \in \pi^{\lambda_{j}}_{j}$, then $\Houtc$ acts on
%  $\mathcal{F}$ by
%\begin{align} \label{HoutAction}
%	f \cdot (A_{1} \otimes \cdots \otimes A_{m}) = & \sum^{r}_{i =
%		1} A_{1} \otimes \cdots \otimes (\Res_{p_{i}} f
%		\varpi_{i}) \cdot A_{i} \otimes \cdots \otimes B_{s}
%		\\ &+ \sum^{s}_{j=1, k=0,1} A_{1} \otimes \cdots
%		\otimes (\Res_{\sigma^{k}_{C}(q_{j})} f \varpi_{j})
%		\cdot B_{j} \otimes \cdots \otimes B_{s}.
%\end{align}
\begin{equation} \label{HoutAction}
f \cdot (A_{1} \otimes \cdots \otimes A_{m}) = \sum^{m}_{i=1} \sum_{p
  \in \nu^{-1}(x_{i})} A_{1} \otimes \cdots \otimes (\Res_{p} f
  \varpi_{i}) A_{i} \otimes \cdots \otimes A_{m} 
\end{equation}
We will say that a meromorphic function $f$ on $C$ is \emph{even} if
$\sigma_{C}^{*}(f) = f$ and \emph{odd} if $\sigma_{C}^{*}(f) = -f$.
If $f$ is even, then $\sum_{p \in \nu^{-1}(x_{i})} \Res_{p} f
\varpi_{i} =0 $, so only odd functions act non-trivially.  Let
$\Houtc^{o}$ denote the space of odd elements in $\mathbb{C}[\Cnot
]$.

\subsection{Coinvariants and conformal blocks} \label{CCb}

Now we give an simpler, alternative definition of the spaces of coinvariants
and conformal blocks for (twisted) $\pi$--modules, extending
Definition 8.1.7 in \cite{FBZ}. 

%\begin{definition}	\label{HeisConfB}
%	The space of \emph{coinvariants} associated to $(C,\{ p_{i}
%\}, \{ q_{j} \}, \{ \pi^{\sigma}_{i} \}, \{ \pi^{\lambda_{j}}_{j} \})$
%is the vector space
%\[
%	\Coinv{r}{s} = \mathcal{F} / \Houtc^{o} \cdot \mathcal{F}.
%\] 
%
%	The space of \emph{conformal blocks}  associated to $(C,\{
%p_{i} \}, \{ q_{j} \}, \{ \pi^{\sigma}_{i} \}, 
%\{ \pi^{\lambda_{j}}_{j} \})$ is the vector space
%\[
%	\Cb{r}{s} = \Hom_{\Houtc^{o}}(\mathcal{F}, \mathbb{C})
%\]
%of $\Houtc^{o}$--invariant functionals on $\mathcal{F}$. 
%\end{definition} 

\begin{definition} \label{HeisConfB}
  The space of \emph{coinvariants} associated to $(X,\{ x_{i} \},
  \{ \pi_{x_{i}} \})$ is the vector space
\[
\Cv{m}= \mc{F} / \Houtc^{o} \cdot \mc{F}
\]
The space of \emph{conformal blocks} associated to $(X,\{ x_{i} \},
 \{ \pi_{x_{i}} \})$ is the vector space
\[
\Cb{m} =  \Hom_{\Houtc^{o}}(\mathcal{F}, \mathbb{C})
\]
of $\Houtc^{o}$--invariant functionals on $\mc{F}$.
\end{definition}

\medskip

\noindent {\bf{Remark.}} The discussion at the end of Section
\ref{Houtc} implies that these definitions remain the same
if we replace $\Houtc^{o}$ by $\Houtc$.\qed

\medskip

We can ask how the space of conformal blocks changes under the
addition of points. Suppose then that to our collection of points
$\points{m}$ we add $x_{m+1}$. Since $\points{m}$ contains all the
branch points of $\nu$, $\nu^{-1}(x_{m+1})$ consists of two points.
Set $\pi_{x_{m+1}} = \pi^{0,p}_{x_{m+1}}$, $p \in
\nu^{-1}(x_{m+1})$. Observe that in the case of the vacuum
representation, $\pi^{0,p}_{x_{m+1}} =
\pi^{0,\sigma_{C}(p)}_{x_{m+1}}$, so there is no choice of point in
the fiber. There exists a natural map:
\begin{equation} \label{CBrestrict1}
\Cb{m+1} \rightarrow \Cb{m}
\end{equation}
given by
\begin{equation} \label{CBrestrict2}
	\phi \rightarrow \phi \vert_{ \bigotimes^{m}_{i=1}
	\pi_{x_{i}} \otimes \vac},
\end{equation}
i.e., we restrict the functional $\phi$ to the vacuum vector in
$\pi_{x_{m+1}}$.  The following lemma, analogous to Proposition 8.3.2 in
\cite{FBZ}, will be proved in Section \ref{CBrProof}.

\begin{lemma}	\label{CBrLemma}
	The map \eqref{CBrestrict2} is an isomorphism.
\end{lemma}

\subsection{Equivalence of Definitions \ref{ConfB2} and \ref{HeisConfB}}

We now have two seemingly different definitions of the space of
conformal blocks: the general Definition \ref{ConfB2} and Definition
\ref{HeisConfB}, which is specific to the case $V=\pi$. We will show
that these two definitions agree.  Applying the Strong Residue Theorem
\ref{ResThm} to our collection of $\End(\mathcal{F})$--valued
one-forms, we obtain

\begin{corollary}
A functional $\phi$ is a conformal block if and only if $\forall
A_{i} \in \pi_{x_{i}}$, the
one-forms
\begin{equation} \label{OneForms}
	\langle \phi, A_{1} \otimes \cdots \otimes ( \varpi_{i} \cdot
	A_{i} ) \otimes \cdots \otimes A_{m} \rangle \in
	\Gamma \left(\coprod_{p \in \nu^{-1}(x_{i}) } \fpd{p},
	\Omega_{\ol{C}} \right) 
\end{equation}
can be extended to a single one-form $\varpi_{\phi}$ on $\Cnot$. 
\end{corollary}

{\bf{Remark.}} $\sigma_{C}$ acts on the space of holomorphic
one-forms. It is clear that $\varpi_{\phi}$ is odd under this
action.\qed

\medskip

We are now ready to prove the equivalence of the two definitions of
conformal blocks. The proof is a generalization of the proof of
Theorem 8.3.3 of \cite{FBZ}.

\begin{theorem} \label{CBEquivalence}
	Let $\phi$ be a linear functional on $\mathcal{F}$ such that
        $\forall A_{i} \in \pi_{x_{i}}$, the one-forms
        \eqref{OneForms} can be extended to a
        single, odd, regular one--form $\varpi_{\phi}$ on
        $\Cnot$. Then the sections
\begin{equation}
	\langle \phi, A_{1} \otimes \cdots \otimes (
	\mc{Y}^{\pi_{x_{i}}} \cdot A_{i} ) \otimes \cdots \otimes
	A_{m} \rangle \in \Gamma( \coprod_{p \in \nu^{-1}(x_{i})}
	\fpd{p},\Pi^{*}_{C} )
	\label{twists2}
\end{equation}
can be extended to a single, invariant, regular section of
$\Pi^{*}_{C}$ on $\Cnot$, and vice versa.
\end{theorem}

\noindent {\em Proof.} From the discussion following Theorem
\ref{PrimaryField}, there exists a map
\[
	\Pi^{*}_{\ol{C}} \rightarrow \Omega_{\ol{C}} 
\]
such that the one-forms \eqref{OneForms} are the projections of
the sections \eqref{twists} on $\discs{i}$ Thus if the
\eqref{twists} extend to $\Cnot$, so will \eqref{OneForms}.
Denote $\Cnot$ by $\Caff$.

Let $\Csq$ denote $\Caff^{2} \backslash \Xi$, where $\Xi$ is the
divisor in $\Caff^{2}$ consisting of pairs $(x,y) \in \Caff^{2}$
satisfying $\sigma^{k}_{C}(x) = \sigma^{l}_{C}(y)$ for some integers
$k,l$. Let $r_{1} : \Csq \rightarrow \Caff$ denote the projection on
the first factor, whose fiber over $q \in \Caff$ is
\[
	\Caff \backslash \Orb{q} 
\]
where $\Orb{q}$ denotes the $H$--orbit of $q$. Let
\[
	\widetilde{\mathcal{O}} = (r_{1})_{*} \mathcal{O}_{\Csq}.
\]
This is a quasi-coherent sheaf on $\Caff$ whose fiber at $q \in \Caff
$ is $\mathbb{C}[ \Caff \backslash \Orb{q} ]$.  Let $\mathcal{G}$
denote the sheaf $\mathcal{F} \otimes \Pi_{C} \vert_{\Caff}$ on
$\Caff$, where we treat $\mathcal{F}$ as a constant sheaf. For $q \in
\Caff$, the fiber $\mathcal{G}_{q}$ is $\mathcal{F} \otimes (\Ptx)_{q}
\cong \mathcal{F} \otimes \pi_{\nu(q)}$, where by $\pi_{\nu(q)}$ we
mean the module along $\nu^{-1}(\nu(q))$ constructed out of the vacuum
$\pi$. The sheaf $\Otil$ acts on $\mathcal{G}$ in such a way that
fibrewise we obtain the action \eqref{HoutAction} of $\mathbb{C}[\Caff
\backslash \Orb{q}] $ on
\[
	\mathcal{F} \otimes \Pi_{\nu(q)} = \bigotimes^{m}_{i=1}
		\pi_{x_{i}} \otimes \pi_{\nu(q)}
\]
Introduce the sheaf of homomorphisms 
\[
	\Ctil = \Hom_{\Otil}(\mathcal{G}, \mathcal{O}_{\Caff})
\]
whose fiber at $q \in \Caff$ is the space of conformal blocks
\[
	\Cbq{m}. 
\]
Lemma \ref{CBrLemma} identifies the fibers of $\Ctil$ with
$\Cb{r}$, thus providing a canonical trivialization of
$\Ctil$. Thus, given $$\phi \in \Cb{m}$$ we obtain a section
$\phitil$ of $\Ctil$ on $\Caff$. Let $\widetilde{\phi}_{q}$ denote
the value of this section at $q \in C$. Let $A_{i} \in \pi_{x_{i}}$ and
\[
	v = \otimes^{m}_{i=1} A_{i} \in \mc{F} 
\]
Contracting $\phitil$ with $v \in \mathcal{F}$, we obtain a section
$\widetilde{\phi}^{v}$ of $\Ptx^{*}$. By construction, the value of
$\phitil^{v}$ on $B \in (\Ptx)_{q} = \pi_{\nu(q)} $ is $\phitil_{q}(v
\otimes B)$.  We wish to show that the section $\phitil^{v}$ is the
analytic continuation of the sections \eqref{twists}.  We begin with
the following lemma.

\begin{lemma} \label{FormLemma}
There exists a regular one-form on $\Caff
\backslash \Orb{q}$ whose restriction to the union $\discs{i}$ is
equal to
\[
	\langle \phitil_{q}, A_{1} \otimes \cdots \otimes \varpi_{i}
		\cdot A_{i} \otimes \cdots \otimes A_{m} \otimes B
		\rangle,
\]
and its restriction to $\coprod_{p \in \nu^{-1}(\nu(q))} \fpd{p}$ is
\[
	\langle \phitil_{q}, A_{1} \otimes \cdots \otimes A_{m}
	\otimes \varpi_{q} \cdot B \rangle
\]
\end{lemma}
\begin{proof}

Since $\phitil_{q}$ is a conformal block, we obtain, using the
definition of the action of $\mathbb{C}[\Caff \backslash \Orb{q}]$ on
$\mathcal{F} \otimes \pi_{\nu(q)}$,
\begin{align*}
	0 = \langle \phitil_{q}, f \cdot (v \otimes B) \rangle = &
	   \sum_{i=1 \cdots m} \sum_{p \in \nu^{-1}(x_{i})} \Res_{p}
	   \langle \phitil_{q}, A_{1} \otimes \cdots \otimes (f
	   \varpi_{i}) \cdot A_{i} \otimes \cdots \otimes B \rangle \\
	   & + \sum_{p=(q,\sigma_{C}(q))} \Res_{p} \langle \phitil_{q},
	   A_{1} \otimes \cdots \otimes A_{m} \otimes \cdots \otimes
	   (f \varpi_{q}) \cdot B \rangle.
\end{align*}
Thus, the strong residue theorem implies that there exists a regular
one-form on $\Caff \backslash \Orb{q}$ having the desired properties.
\end{proof}

\begin{remark}
  If $x \in X, p \in \nu^{-1}(x)$, $\nu^{-1}(x)$ consists of
  $\frac{2}{N}$ points, and $z^{\frac{1}{N}}$ is a formal
  coordinate at $p$, then the endomorphism-valued one-form $\varpi$ on
  $\fpd{p}$ has the expression
\[
\varpi = N z^{\frac{N-1}{N}} b(\zf{N}) d \zf{N}
\]
where
\[
b(\zf{N}) = \sum_{n \in \frac{1}{N} + \Z} b_{n} z^{-n - 1}
\]
\end{remark}

\noindent Now we prove Theorem \ref{CBEquivalence}.  Suppose that $q$
is near $p \in \nu^{-1}(x_{i})$.  Choose a small analytic neighborhood
$U$ of $p$ with special coordinate $\zf{N}_{i}$ centered on $p$, such
that $q \in U_{i}$. If $\nu$ is unramified, $N=1$, and $\zf{N}_{i}$ is
any coordinate centered at $p$, otherwise $N=2$ and $\zf{2}_{i}$ is a
$\sigma_{C}$--special coordinate. Since $q \ne p$, $w=z_{i}-z_{i}(q)$
is a coordinate centered at $q$, in some neighborhood $W$ of $q$. Near
$q$, we can trivialize
\[
	\pi \times W \cong \Ptx \vert_{W}  
\]
via
\[
	(B, q) \rightarrow (B,z_{i}-z_{i}(q)).
\]
It remains to prove the following: 

\begin{lemma} 
	$\forall B \in \pi$, \label{AnalContLemma}
\begin{align*}
	\langle \phitil_{q}, v \otimes B \rangle & = \langle
		\phitil_{q}, A_{1} \otimes \cdots \otimes
		Y^{\pi_{x_{i}}}(B, q) \cdot A_{i} \otimes \cdots
		\otimes A_{m} \otimes \vac \rangle \\ & = \langle
		\phi, A_{1} \otimes \cdots \otimes
		Y^{\pi_{x_{i}}} (B, q) \cdot A_{i} \otimes \cdots
		\otimes A_{m} \rangle.
\end{align*}
\end{lemma} 

\begin{proof} 
	The second equality follows from Lemma \ref{CBrLemma}. The
first equality is proved by induction. It obviously holds for $B =
\vac$. Now, denote by $\pi^{(r)}$ the subspace of $\pi$ spanned by all
monomials of the form $\widetilde{b_{i_{1}}} \cdots
\widetilde{b}_{i_{k}} \vac$, where $k \leq r$. Suppose that we have
proved Lemma \ref{AnalContLemma} for all $B \in \pi^{(r)}$. The
inductive step is to prove the equality for elements of the form $B' =
\widetilde{b_{n}} \cdot B$. By our inductive hypothesis, we know that
if $B \in \pi^{(r)}$, then

\begin{multline*}
	\langle \phitil_{q}, A_{1} \otimes \cdots \otimes ( N
		z^{\frac{N-1}{N}}_{i} b(\zf{N}_{i})) \cdot A_{i}
		\otimes \cdots \otimes B \rangle d \zf{N}_{i} \\ =
		\langle \phi, A_{1} \otimes \cdots \otimes N  
		Y^{\pi_{x_{i}}}(B,q) z^{\frac{N-1}{N}}_{i} b(\zf{N}_{i})
		\cdot A_{i} \otimes \cdots \otimes A_{m} \rangle d
		\zf{N}_{i}.
\end{multline*}
According to Lemma \ref{FormLemma} above, we also have
\[
	\langle \phitil_{q}, A_{1} \otimes \cdots \otimes ( N
		z^{\frac{N-1}{N}}_{i} b(\zf{N}_{i})) \cdot A_{i}
		\otimes \cdots \otimes A_{m} \otimes B \rangle d
		\zf{N}_{i} = \langle \phitil_{q}, A_{1} \otimes \cdots
		\otimes A_{m} \otimes \widetilde{b}(w) \cdot B \rangle
		dw.
\]
Using locality and associativity, we obtain
\begin{align*}
	\langle \phi, A_{1} \otimes \cdots & \otimes
		Y^{\pi_{x_{i}}}(C,q) b(\zf{N}_{i}) \cdot A_{i}
		\otimes \cdots \otimes A_{m} \rangle  N
		z^{\frac{N-1}{N}}_{i} d \zf{N}_{i} \\ &=\langle \phi,
		A_{1} \otimes \cdots \otimes b(\zf{N}_{i})
		Y^{\pi_{x_{i}}}(C,q) \cdot A_{i} \otimes \cdots
		\otimes A_{m} \rangle N z^{\frac{N-1}{N}}_{I} d \zf{N}_{i}
		\\ &=\langle \phi, A_{1} \otimes \cdots \otimes
		Y^{\pi_{x_{i}}}(\widetilde{b}(z_{i} - z_{i}(q))
		\cdot B, q) \cdot A_{i} \otimes \cdots \otimes
		A_{m}\rangle N z^{\frac{N-1}{N}}_{i} d \zf{N}_{i} \\ &=
		\langle \phi, A_{1} \otimes \cdots \otimes
		Y^{\pi_{x_{i}}}(\widetilde{b}(w) \cdot B, q )
		\cdot A_{i} \otimes \cdots \otimes A_{m} \rangle N
		z^{\frac{N-1}{N}}_{i} d \zf{N}_{i}\\ &= \langle \phi,
		A_{1} \otimes \cdots \otimes
		Y^{\pi_{x_{i}}}(\widetilde{b}(w) \cdot B, q )
		\cdot A_{i} \otimes \cdots \otimes A_{m} \rangle dw,
\end{align*}
where the last step holds because $ N z^{\frac{N-1}{N}}_{i} d \zf{N}_{i} =
		dw$. Combining these relations, we obtain
\[
	\langle \phitil_{q}, A_{1} \otimes \cdots \otimes A_{m}
		\otimes \widetilde{b}(w) \cdot B \rangle dw = \langle
		\phi, A_{1} \otimes \cdots \otimes
		Y^{\pi_{x_{i}}}(\widetilde{b}(w) \cdot B, q)
		\cdot A_{i} \otimes \cdots \otimes A_{m} \rangle dw.
\]
Multiplying both sides by $w^{n}$ and taking residues, we find that
\[
	\langle \phitil_{q}, A_{1} \otimes \cdots \otimes A_{m}
		\otimes \widetilde{b}_{n} \cdot B \rangle = \langle
		\phi, A_{1} \otimes \cdots \otimes
		Y^{\pi_{x_{i}}}(\widetilde{b}_{n} \cdot B, q)
		\cdot A_{i} \otimes \cdots \otimes A_{m} \rangle.
\]
Equivariance of the sections \ref{twists2} follows from the fact that they
are invariant on all $\fpd{p}$ where $p$ is a fixed point of
$\sigma_{C}$.  This completes the proof of Theorem
\ref{CBEquivalence}.

\end{proof}

\subsection{Proof of Lemma \ref{CBrLemma}} \label {CBrProof}

We start with the following fact.

\begin{lemma} \label{PrincPartLemma}

Let $p \in \nu^{-1}(x_{m+1})$. For every principal part $f_{-}$ at
$p$, there exists an odd function $f \in \mathbb{C}[\Caff \backslash
\Orb{p} ]$ whose principal part at $p$ is $f_{-}$.

\end{lemma}

\begin{proof}

Let $D$ be an effective divisor symmetric under the action of
$\sigma_{C}$ (i.e., if $D = \Sigma c_{q} \cdot q$, then $c_{q} =
c_{\sigma_{C}(q)}$ ), supported on $\{ \nu^{-1}(x_{i}) \}_{i=1 \cdots
m}$.  Denote the canonical divisor of $C$ by $K_{C}$. For $\deg(D) >
\deg(K_{C})$, the Riemann-Roch theorem implies that
\[
	\on{dim} \mathcal{L}(D) = \deg(D) + 1 - g_{C}.
\]
It follows that 
\[
	Q_{n+1} = \mathcal{L}(D+ (n+1) \cdot p + (n+1)
	\sigma_{C}(p)) / \mathcal{L}(D + n \cdot p + n
	\cdot \sigma_{C}(p))
\]
for $n \geq 0$ is two-dimensional. Furthermore, $Q_{n+1}$ carries an
action of $\sigma_{C}$. Suppose now that $Q_{n+1}$ is spanned by the
images of two even functions $f_{1}, f_{2}$. Since $f_{i}$ are even,
they have poles of the same order at both $p$ and
$\sigma_{C}(p)$, and so will any linear combination. But this
contradicts the fact that
\[
	\mathcal{L}(D + (n+1) p + n \sigma_{C}(p) ) /
	\mathcal{L}(D + n \cdot p + n \cdot \sigma_{C}(p))
\]
is one-dimensional. It follows that for each $n \geq 0$, $Q_{n+1}$
contains an odd function.
\end{proof}

Now we prove the statement equivalent to Lemma \ref{CBrLemma} that the
corresponding spaces of coinvariants 
\[
\Hco{m} \quad \on{and} \quad \Hco{m+1}
\]
are isomorphic. Recall that $\pi_{x_{m+1}}$ here is the vacuum module
$\pi^{0,p}_{x_{m+1}}$. Let $\Caffp = \Caff \backslash \Orb{p}$, and
$\Houtcp = \mathbb{C}[\Caff' ]$. The space
\[
\Hco{m} \qquad (\on{resp.} \Hco{m+1}
\]
is identified with
the $0$th homology of the Lie algebra $\Houtc^{o}$
(resp. $\Houtcp^{o}$) with coefficients in $\mathcal{F}$
(resp. $\mathcal{F} \otimes \pi_{x_{m+1}}$). Lemma \ref{PrincPartLemma}
implies that the sequence
\begin{equation} \label{HoutSequence}
	0 \rightarrow \Houtc^{o} \rightarrow \Houtcp^{o}
	\xrightarrow{\mu} w^{-1}_{s+1} \mathbb{C}[w^{-1}_{s+1} ]
	\rightarrow 0
\end{equation}
is exact, where $\mu$ is the map that attaches to a function its
principal part at $p $. The homology of $\Houtcp^{o}$ with
coefficients in $\mathcal{F} \otimes \pi_{x_{m+1}}$ is computed using the
Chevalley complex
\[
	C^{\bullet} = \mathcal{F} \otimes \pi_{x_{m+1}} \otimes
	\bigwedge{}^{\bullet} (\Houtcp^{o})
\]
with the differential $d: C^{i} \rightarrow C^{i-1}$ given by the formula
\[
	d = \sum_i f_{i} \otimes \psi^{*}_{i}
\]
where $\{ f_{i} \}$ is a basis in $\Houtcp^{o}$ and $\{ \psi^{*} \}$
is the dual basis of $(\Houtcp^{o})^{*}$ acting on $
\bigwedge^{\bullet} (\Houtcp^{o})$ by contraction.

Choose pull-backs $z_{n}, n < 0$, of $w^{n}_{s+1}, n < 0$, in
$\Houtcp^{o}$ under $\mu$. Because of the exactness of the sequence
\eqref{HoutSequence}, we can choose a basis $\{ f_{i} \}$ in
$\Houtcp^{o}$ which is a union of $\{ z_{n} \}_{n < 0}$, and a basis
of $\Houtc^{o}$. In this basis we may decompose
\[
	d = d_{\Caff} + \sum z_{n} \otimes \phi^{*}_{n},
\]
where $d_{\Caff}$ is the differential for $\Houtc^{o}$, and
$\phi^{*}_{n}$ denotes the element of the dual basis to $\{ f_{i} \}$
corresponding to $z_{n}$.

We need to show that the homologies of this complex are isomorphic to
the homologies of the complex $\mathcal{F} \otimes
\bigwedge^{\bullet}(\Houtc^{o})$.  Introduce an increasing filtration
on $\pi_{x_{m+1}}$, letting $\pi^{(r)}_{m+1}$ be the span of all
monomials of order less than or equal to $m$ in $\wt{b}_{n}, n <
0$. Now introduce a filtration $\{ F_{i} \}$ on the Chevalley complex
$C^{\bullet}$ by setting
\[
	F_{i} = \on{span}\{ v \otimes B \otimes D \vert v \in
	\mathcal{F}, B \in \pi^{(m)}_{x_{m+1}}, D \in
	\bigwedge{}^{i-m}(\Houtcp^{o}) \}.
\]
Our differential preserves this filtration.

Consider now the spectral sequence associated to the filtered complex
$C^{\bullet}$. The zeroth term $E^{0}$ is the associated graded space
of the Chevalley complex, isomorphic to
\[
	(\pi_{s+1} \otimes \bigwedge{}^{\bullet} (\phi^{*}_{n})_{n <
	0} ) \otimes (\mathcal{F} \otimes
	\bigwedge{}^{\bullet}(\Houtc^{o})).
\]
The zeroth differential acts along the first factor of the above
decomposition, and is given by the formula
\[
	d^{0} = \sum_{n < 0} b_{n} \otimes \phi^{*}_{n},
\]
because on the graded module the operator $z_{n}$ acts as $b_{n}, n <
0$. But $\pi_{s+1}$ is isomorphic to the symmetric algebra with
generators $b_{n}, n < 0$, and our differential is simply the Koszul
differential for this symmetric algebra. It is well-known that the
zeroth homology of this complex is isomorphic to $\mathbb{C}$, and all
other homologies vanish. Therefore all positive homologies of $d^{0}$
vanish, while the zeroth homology is $ \mathcal{F} \otimes
\bigwedge^{\bullet}(\Houtc) $. Hence, the $E^{1}$ term coincides as a
vector space with the Chevalley complex of the homology of
$\Houtc^{o}$ with coefficients in $\mathcal{F}$. Also, the $E^{1}$
differential coincides with $d_{\Caff}$, which is the corresponding
Chevalley differential. We thus obtain the desired isomorphism
\[
	H_{i}(\Houtcp^{o}, \mathcal{F} \otimes \pi_{m+1}) \cong
	H_{i}(\Houtc^{o}, \mathcal{F}).
\]

\section{Affine vertex algebras}

In Section \ref{CCb} we have shown that in the case of the Heisenberg
vertex algebra the space of conformal blocks had a simple realization
as the dual of a certain space of twisted coinvariants. In this
section we present a similar realization in the case of vertex
algebras attached to affine Kac-Moody algebras.

\subsection{The vacuum module $\Vkg$}

Let $\g$ denote a complex simple Lie algebra, $\Lg = \g \otimes [t,
t^{-1}]$ its loop algebra, and $\ghat$ the corresponding affine
Kac-Moody Lie algebra. For $k \in \mathbb{C}$, let $\mathbb{C}_{k}$
denote the one-dimensional representation of $\g[t] \oplus \mathbb{C}
\cdot K$ where $\g[t]$ acts by 0, and $K$ acts by $k$.  It is well
known that the \emph{vacuum module}
\[
	\Vkg = \Ind^{\ghat}_{\g[t] + \mathbb{C} \cdot K}
	\mathbb{C}_{k}
\]
has the structure of a vertex algebra (see for instance Section 3.4.2
of \cite{FBZ}). 

Pick a basis $\{ J^{a} \}_{a = 1 \cdots d}$ (where $d = \dim(\g)$) of
$\g$, and let $\{ J_{a} \}_{a=1 \cdots d}$ be its dual basis with
respect to the normalized Killing form. Suppose that $k \neq -\dcox$
(where $\dcox$ is the dual Coxeter number of $\g$) and set
\[
	S = \frac{1}{2(k + \dcox)} \sum^{d}_{a=1} ( J_{a} \otimes
	t^{-1}) ( J^{a} \otimes t^{-1}) \vac.
\]
This is the \emph{Sugawara vector} which determines a conformal
structure on $\Vkg$ when $k \ne - \dcox$. In what follows, we will
always use this conformal structure on $\Vkg$.

Let $\sigma_{\g}$ be an automorphism of $\g$ of finite order $N$. Then
$\sigma_{g}$ induces a conformal automorphism of $\Vkg$, which we will
denote by $\sigma_{\Vkg}$.

In particular, consider the case when $\sigma_\g$ is an outer
automorphism (note that this is not necessary for the results
below). Thus, $N=2$ when $\g= A_{n}, D_{m}, m \ne 4, E_{6}$, and $N=3$
when $\g = D_{4}$. The following result is proved in \cite{LI}

\begin{lemma}
The $\sigma_{\Vkg}$--twisted $\Vkg$--modules are precisely the
$\ghat^\sigma $--modules from the category $\OO$, where
$\ghat^{\sigma}$ is the twisted affine Kac-Moody algebra associated to
the automorphism $\sigma_{\g}$.
\end{lemma}

\subsection{The Lie algebra $\goutc$}

We keep the notation of Section \ref{HeisCB}.  Let $C$ be an algebraic
curve with an automorphism $\sigma_{C}$ of order $N$
(where $N=2$ or $3$ depending on $\g$), and let
$\{x_{i} \}_{i=1 \cdots m}$ be a collection of points of $X$
containing the branch points of $\nu$.  Denote $\Cnot$ by $\Caff$. Let
us write $\g = \oplus^{N-1}_{l=0} \g_{l}$, where $\g_{l}$ denotes the
eigenspace of $\sigma_{\g}$ corresponding to the eigenvalue
$e^{\frac{2 \pi i l }{N}}$. Then $\sigma_{C}$ acts on
$\mathbb{C}[\Caff]$--the ring of functions on $\Caff$, and so we can
write $ \mathbb{C}[\Caff] = \oplus \mathbb{C}[\Caff]_{l} $, where $
\mathbb{C}[\Caff]_{l}$ consists of those functions $f$ such that
$\sigma^{*}_{C}(f) = e^{\frac{2 \pi i l}{N}} f$. Let
\[
	\goutc = \oplus^{N}_{l=1} \left( \g_{l} \otimes
	\mathbb{C}[\Caff]_{l} \right).
\]

\subsection{Coinvariants and conformal blocks}

For $x \in X$, $V$--modules along $\nu^{-1}(x)$ can be constructed
from ordinary or twisted $V$--modules using the same technique that
was used in the Heisenberg case in Section \ref{Houtc}. More
precisely, if $x$ is a branch point of $\nu$,$p=\nu^{-1}(x)$, and
$\sigma_{C,p}$ is the monodromy around $x$, then any
$\sigma_{C,p}$--twisted $\Vkg$--module gives rise to a $\Vkg$--module
along $\nu^{-1}(x)$. Similarly, if $\nu^{-1}(x)$ consists of $N$
points, then an ordinary $\Vkg$--module and a choice of point $p \in
\nu^{-1}(x)$ gives rise to a $\Vkg$--modules along $\nu^{-1}(x)$.

Let $\{ \mc{M}_{x_{i}}\}$ be a collection of $\Vkg$--modules along $\{
x_{i} \}$ constructed in this manner. Thus for each $x_{i}$, we have a
distinguished point $p_{i} \in \nu^{-1}(x_{i})$. Pick special
coordinates $\zf{N_{i}}$ near $p_{i}$, where $N_{i}=1$ if $\nu$ is
unramified at $p_{i}$ and $N_{i}=N$ otherwise. Set 
\[
	\mathcal{F} = \bigotimes^{r}_{i=m} \mc{M}_{x_{i}}
\]
Then $\goutc$ acts on $\mc{F}$ as follows:
\[
	h \cdot (A_{1}  \otimes \cdots \otimes A_{m}) = \sum_{i} A_{1}
		      \otimes \cdots \otimes [h]_{p_{i}} \cdot A_{i}
		      \otimes \cdots \otimes A_{m} 
\]
where $[h]_{p}$ denotes the Laurent series expansion of $h$ around $p
\in C$ in the special coordinate that was selected.

We are now ready to give an alternative, simplified definition of
twisted coinvariants and conformal blocks for $V_k(\g)$, extending 
the definition of Section 8.2.1 in \cite{FBZ}:

\begin{definition} \label{KMCbDef}	
The space of \emph{coinvariants} is the vector space
\[
	\wt{\mc{H}}_{\Vkg}(X,\{ x_{i} \}, \mc{M}_{x_{i}})_{i=1 \cdots
	m} = \mc{F} / \goutc \cdot \mc{F}
\]
The space of \emph{conformal blocks} is its dual:
\[
	\wt{\mc{C}}_{\Vkg}(X, \{ x_{i} \}, \mc{M}_{x_{i}})_{i=1 \cdots
	m} = \on{Hom}_{\goutc}(\mc{F}, \mathbb{C})
\]
\end{definition} 

\medskip

The following theorem is proved using the same methods as Theorem
\ref{CBEquivalence}. 

\begin{theorem}
In the case of the vertex algebra $\Vkg$, Definition \ref{ConfB2}
is equivalent to Definition \ref{KMCbDef}. 
\end{theorem}

\end{document}